\documentclass[12pt]{amsart}
\usepackage{amssymb}
\usepackage[all]{xy}
\usepackage[toc,page,title,titletoc,header]{appendix}
\usepackage{xcolor}

\theoremstyle{plain}
\newtheorem{thm}{Theorem}[section]

\newtheorem*{thm1.1}{Theorem 1.1}

\newtheorem{lemma}[thm]{Lemma}
\newtheorem{lem}[thm]{Lemma}
\newtheorem{cor}[thm]{Corollary}
\newtheorem{pro}[thm]{Proposition}

\theoremstyle{definition}

\newtheorem{rem}[thm]{Remark}

\newtheorem{defi}[thm]{Definition}
\newtheorem{exe}[thm]{Example}

\newtheorem{que}[thm]{Question}
\newtheorem{con}[thm]{Conjecture}

\newtheorem*{dmlcon}{Dynamical Mordell-Lang Conjecture}

\numberwithin{equation}{section}
\newcounter{elno}                

\newcommand{\la}{\lambda}

\newcommand{\Supp}{{\rm Supp}}
\newcommand{\rank}{{\rm rank}}

\newcommand{\Pic}{{\rm Pic}}

\newcommand{\Char}{{\rm char}}

\newcommand{\id}{{\rm id}}

\newcommand{\Gal}{{\rm Gal}}

\newcommand{\boxtensor}{{\Box\kern-9.03pt\raise1.42pt\hbox{$\times$}}}

\newcommand{\hotimes}{\hat{\otimes}}

\newcommand{\supp}{{\rm supp}\,}

\newcommand{\propsubset}
{\mbox{$\textstyle{
\subseteq_{\kern-5pt\raise-1pt\hbox{\mbox{\tiny{$/$}}}}}$}}

\newcommand{\sM}{{\mathcal M}}

\newcommand{\sP}{{\mathcal P}}

\newcommand{\A}{{\mathbb A}}

\newcommand{\C}{{\mathbb C}}
\newcommand{\D}{{\mathbb D}}

\newcommand{\F}{{\mathbb F}}

\renewcommand{\P}{{\mathbb P}}
\newcommand{\Q}{{\mathbb Q}}
\newcommand{\R}{{\mathbb R}}

\newcommand{\Z}{{\mathbb Z}}
\newcommand{\bk}{{\mathbf{k}}}

\begin{document}
\title[]{Around the dynamical Mordell-Lang conjecture}

\author{Junyi Xie}

\address{Beijing International Center for Mathematical Research, Peking University, Beijing 100871, China}

\email{xiejunyi@bicmr.pku.edu.cn}




\thanks{The author is supported by the NSFC Grant No.12271007.}

\date{\today}

\bibliographystyle{alpha}


\maketitle

\begin{abstract}There are three aims of this note. The first one is to report some advances around the dynamical Mordell-Lang (=DML) conjecture. Second, we generalize some known results. For example, the Dynamical Mordell-lang conjecture was known for endomorphisms of $\A^2$ over $\overline{\Q}$. We generalize this result to all endomorphisms of $\A^2$ over $\C$. We generalize the weak DML theorem to a uniform version and to a version for partial orbit. Using this, we give a new proof of the Kawaguchi-Silverman-Matsuzawa's upper bound for arithmetic degree. We indeed prove a uniform version which works in both number field and function field case in any characteristic and it works for partial orbits. We also reformulate the ``$p$-adic method", in particular the $p$-adic interpolation lemma in language of Berkovich space and get more general statements.
The third aim is to propose some further questions. 
\end{abstract}

\tableofcontents

\section{Introduction}

The dynamical Mordell-Lang conjecture describes the return time set for an algebraic dynamical system. 
It was proposed by Ghioca and Tucker in~\cite{Ghioca2009}, under the influence of Shou-Wu Zhang, see \cite[top of Page 306]{Ghioca2008a}.
Its original version only concerns endomorphisms. We present here a slightly more general form for rational self-maps.


Let $\bk$ be an algebraically closed field.
Let $X$ be a quasi-pro\-jective variety over $\bk$. Let
$f: X \dashrightarrow X$ be a dominant rational self-map. Let $V$ be a Zariski closed subset of $X$. 
Denote by $I(f)$ the inderterminacy locus of $f$. 
Let $X_f(\bk)$ be the set of points in $X(\bk)$ whose $f$-orbit $O_f(x):=\{n\geq 0|\,\, f^n(x)\}$ is well-defined i.e. for every $n\geq 0$, $f^n(x)\not\in I(f)$.

\begin{dmlcon}[=DML]\label{dml}
	If $\Char\,\bk=0$, then for every $x\in X_f(\bk)$  the set $\{n\in \mathbb{N}|\,\,f^n(x)\in V(\bk)\}$ is a union of at most finitely many arithmetic progressions\footnote{An arithmetic progression is a set of the form $\{an+b|\,\,n\in \mathbb{N}\}$ with $a,b\in \mathbb{N}$. In particular, when $a=0$, it contains only one point.}.
\end{dmlcon}

\rem
The assumption that $\bk$ is algebraically closed is just for the convenience. In general, easy to see that DML holds over $\overline{\bk}$ implies the same statement over $\bk$.
Moreover, we may further assume that $\bk=\C$. Indeed, DML over $\C$ implies DML over any field of characteristic $0$. 
\endrem

In practice, the following geometric form of DML is often easier to handle and apply since it explains the mechanism to make the return time set infinite.
\begin{con}[=geometric form of DML]\label{dmlgeo}
	Assume $\Char\,\bk=0$ and $V$ is irreducible of dimension at least $1$. For $x\in X_f(\bk)$ if $O_f(x)\cap V$ is Zariski dense in $V$, then $V$ is periodic i.e. the orbit of $\eta_V$ is well-defined and $f^m(V)\subseteq V$ for some $m\geq 1.$
\end{con}
Here $\eta_V$ is the generic point of $V$ and  $f^m(V):=\overline{f^m(\eta_V)}.$
Not hard to see that DML and its geometric form are equivalent. The DML conjecture fits the general principle of ``unlikely intersection problem" i.e.  
an irreducible subvariety containing a Zariski dense subset of ``special" points should be ``special" \cite{Zannier2012}. In the setting of DML, a point is ``special" if it is contained in the orbit of $x$ and an irreducible subvariety is ``special" if it is periodic.

\medskip

\begin{defi}We say that DML holds for $(X,f)$ if the statement of DML conjecture holds for $(X,f)$.
Further, for every $i=0,\dots, \dim X$, we say that DML($i$) holds for $(X,f)$ if the statement of Conjecture \ref{dmlgeo} holds for every irreducible subvariety of dimension $\leq i.$
\end{defi}
It is clear that DML($0$) holds and DML=DML($\dim X-1$)=DML($\dim X$) for $(X,f)$.

\medskip

\subsection{The plan of the note}
The first aim of this note is to report some advances around the DML conjecture;
the second is to generalize some known results and
the third aim is to propose some further questions. 

\medskip

The note is organized as follows:

In Section \ref{origindml}, we discuss the two origins of DML: Diophantine geometry and recurrence sequences.

In Section \ref{knownresults}, we recall some known results of the DML conjecture. Moreover, we show that how to modify the proof of 
 \cite{Xie2017a}, to generalize it to all endomorphisms $\A_{\C}^2$.
\begin{thm}[=Theorem \ref{thmendoplaneoverc}]\label{thmendoplaneoverci}
The DML conjecture holds when $f$ is an endomorphism of $\A^2_{\C}.$
\end{thm}

In Section \ref{dmlpositivec}, we recall some advances of DML in positive characteristic.

In Section \ref{weakdml}, we discuss the Weak DML problem which weaken the DML conjecture in time.
We first recall the Weak DML theorem of Bell, Ghioca and Tucker \cite[Corollary 1.5]{Bell2015}.
Then we generalize this result to a uniform version. 
\begin{thm}[Uniform Weak DML=Theorem \ref{thmdmlweaku}]\label{thmdmlweakui}
Let $X$ be a variety over a field $\bk$. Let $f: X\dashrightarrow X$ be a rational self-map. 
Let $V$ be a proper closed subset of $X$. 
For every $\epsilon>0$, there is $N\geq 1$ such that for every $x\in X_f(\bk)$ with Zariski dense orbit and every interval $I$ in $\Z_{\geq 0}$ with $|I|\geq N,$
we have $$\frac{\#(\{n\geq 0|\,\, f^n(x)\in V\}\cap I)}{\#I}<\epsilon.$$
\end{thm}
We also prove a version of Weak DML for partial orbit (c.f. Theorem \ref{thmperiodicapp}).

In Section \ref{sectionksmuupb}, we apply the Weak DML (for partial orbit) to prove a uniform version of the Kawaguchi-Silverman-Matsuzawa's upper bound.
\begin{thm}[Uniform KSM's upper bound=Theorem \ref{thmuniKSM'sineq}]\label{thmuniKSM'sineqi}
Let $\bk$ be either $\overline{\Q}$ or the algebraic closure of a function field $K(B)$ for a smooth projective curve $B$ over an algebraically closed field $K.$
Let $f: X\dashrightarrow X$ be a dominant rational map defined over $\bk$. Let $h$ be any Weil height on $X$ associated to some ample line bundle.
We denote by $h^+:=\max\{h,1\}.$

Then for
any $\epsilon>0$, there exists $C>0$ such that
$$h^+(f^n(x))\leq C(\la_1(f)+\epsilon)^nh^+(x)$$
for all $n\geq 0$ and $x\in X_f(\bk, n)$. In particular, for any $x\in X_f(\bk)$, we have $$\overline{\alpha_f}(x)\leq \la_1(f),$$
where $\overline{\alpha_f}(x)$ stands for the upper arithmetic degree.
\end{thm}

In Section \ref{almostdml}, we discuss the almost DML, which weaken the DML in space.

In Section \ref{sectionautaff}, we reformulate the ``$p$-adic method", in particular the $p$-adic interpolation lemma using the language of Berkovich spaces. We also generalize the results in more general setting.

In Section \ref{question}, we list some questions relate to the DML conjecture.


\subsection*{Acknowledgement}
The author would like to thank Yohsuke Matsuzawa who asked the question whether the KSM's upper bound of arithmetic degree on singular varieties can be promoted to a uniform version. This motivates the current form of Theorem \ref{thmuniKSM'sineq}.
The author would like to thank Ariyan Javanpeykar for his help on Siegel's theorem over finitely generated domain. 
The author would like to thank Thomas Tucker and She Yang for their helpful comments on the first version of this work.
The author would like to thank the Simons foundation and the organizers, Laura DeMarco and Mattias Jonsson, for inverting me to the Simons symposium ``Algebraic, Complex and arithmetic Dynamics (2019)" and for providing the chance to write this note.

\section{Origins of DML}\label{origindml}
\subsection{Diophantine geometry}
A motivation of the DML conjecture is the
Mordell-Lang conjecture on subvarieties of semiabelian varieties which is now a theorem of Faltings \cite{Faltings1994} and Vojta \cite{Vojta1996}.
\begin{thm}[=Mordell-Lang conjecture]
	Let $V$ be a subvariety of a semiabelian variety $G$ over $\C$ and let $\Gamma$ be a finitely generated subgroup of $G(\C)$. Then $V(\C)\bigcap \Gamma$ is a union of at most finitely many translates of
	subgroups of $\Gamma$.
\end{thm}
In fact, according to the well-known dictionary between arithmetic dynamics and Diophantine geometry,
the DML conjecture is an analogy of the classical Mordell-Lang conjecture, see \cite[beginning of Preface]{Silverman2007} and \cite[\S 6.5]{Silverman2012}.
Easy to see that the DML conjecture for translations on semiabelian varieties implies the classical Mordell-Lang conjecture in the case when $\Gamma$ has rank $\leq 1$.

\medskip

It will be interesting to have a nice formulation of a higher rank analogy of the DML conjecture. One note that the naive formulation should not work as shown in the following example:
\begin{exe}
Let $f,g$ be automorphisms of $\A^2_{\C}$ defined by $f:(x,y)\mapsto (x+1,y)$ and $g:(x,y)\mapsto (x,y+1)$.
One note that $f$ and $g$ commute to each other.
Let $V$ be the curve $\{y=x^2\}$ and $p:=(0,0)$. Then the set $\{(m,n)\in \Z_{\geq 0}^2|\,\, f^m\circ g^n(p)\in V\}$ is infinite, but contains no
infinite translation of subsemi-group of $\Z_{\geq 0}^2.$
\end{exe}
The classical Mordell-Lang conjecture  has the same
phenomenon for additional groups.
It does not hold 
when $\rank\, \Gamma\geq 2$, but holds when $\rank\, \Gamma\leq 1.$

\medskip

\subsection{DML and recurrence sequences}
Another origin of the DML conjecture is the Skolem-Mahler-Lech Theorem \cite{Lech1953} on linear recurrence sequences.
\begin{thm}[Skolem-Mahler-Lech Theorem]\label{thmsml}
Let $\{A_n\}_{n\geq 0}$ be any recurrence sequence satisfying $$A_{n+l}=\sum_{i=0}^{l-1}a_{i}A_{n+i}$$ for all $n\geq 0$, where $l\geq 1$. Then the set $\{n\geq 0|\,\,A_n=0\}$ is a union of at most finitely many arithmetic progressions.
\end{thm}

This statement is equivalent to the DML conjecture for the linear map $$f:(x_0,\dots,x_{l-1})\mapsto(x_1,\dots,x_{l-1},\sum_{i=0}^{l-1}a_{i}x_i)$$
and for the hyperplane $V=\{x_0=0\}.$

\bigskip

There are two natural ways to generalize the Skolem-Mahler-Lech Theorem as follows, both of them are subsequences of the DML conjecture.

\medskip

First, if we allow the recurrence relation in Theorem \ref{thmsml}  to be any polynomial, we get a non-linear generalization of the Skolem-Mahler-Lech Theorem.
\begin{con}[=non-linear SML]\label{connonlinearsml}
Let $\{A_n\}_{n\geq 0}$ be any recurrence sequence satisfying $$A_{n+l}=F(A_n,\dots, A_{n+l-1})$$ for all $n\geq 0$, where $l\geq 1$ and $F\in \C[x_0,\dots, x_{l-1}]$. Then the set $\{n\geq 0|\,\,A_n=0\}$ is a union of at most finitely many arithmetic progressions.
\end{con}
Conjecture \ref{connonlinearsml} is implied by the DML conjecture for the polynomial map $$f:(x_0,\dots,x_{l-1})\mapsto(x_1,\dots,x_{l-1},F(x_0,\dots,x_l))$$
where $l\geq 1$, $F\in \C[x_0,\dots, x_{l-1}]$
and for the hyperplane $V=\{x_0=0\}.$

\medskip

By \cite[Theorem 1]{Xie2017a}, the non-linear SML holds when $l= 2$ and $F$ is defined over $\overline{\Q}.$
We will see later in Theorem \ref{thmendoplaneoverc}, 
with a few modification of the proof of \cite[Theorem 1]{Xie2017a}, it can be generalized to 
any endomorphism of $\A_{\C}^2$.
Hence the non-linear SML holds when $l= 2$.

\medskip

Another way to generalize Theorem \ref{thmsml} is allowing the coefficients $a_i$ to be non-constant. 
If we ask $a_i$ to be generated by iterating of a rational function,  we get the following conjecture.
\begin{con}[=SML with non-constant coefficients]\label{connonconsml}
Let $l\geq 1$ and let $g, a_i\in \C(x), i=0,\dots, l-1$ be rational functions.
Let $\alpha\in \C$ be a point such that 
$a_i(g^n(\alpha))\neq \infty$ for every $i=0,\dots, l-1, n\geq 0.$ Let $\{A_n\}_{n\geq 0}$ be any recurrence sequence satisfying $$A_{n+l}=\sum_{i=0}^{l-1}a_{i}(g^n(\alpha))A_{n+i}$$ for all $n\geq 0$. Then the set $\{n\geq 0|\,\,A_n=0\}$ is a union of at most finitely many arithmetic progressions.
\end{con}

\rem
In the case where $g(x)=x+1$, such a sequence $A_i$ is called a \emph{p-recursive} sequence.  
They are exactly the coefficients of $D$-finite formal power series.
Such sequences play a very important role in symbolic computation \cite{Stanley1980}.
In this case, Conjecture \ref{connonconsml} becomes Rubel's conjecture \cite[Question 16]{Rubel1983}.
\endrem

As shown in \cite[Page 4]{Ghioca2018}, Conjecture \ref{connonconsml} is implied by the DML conjecture for certain skew linear map $f: \P^1\times \A^l\dashrightarrow \P^1\times \A^l$.
By \cite[Theorem 1.2]{Ghioca2018},  Conjecture \ref{connonconsml} holds when $\deg g\geq 2.$
When $\deg g=1$, see \cite[Theorem 1.2]{Ghioca2018} for a weaker result.

\medskip

Conjecture \ref{connonconsml} strongly relates to the Picard-Vessiot problem in difference Galois theory. 
Keep the notations in Conjecture \ref{connonconsml}.
 Then $g$ defines a
difference field $(\bk(x),\sigma)$ where sigma is the endomomorphism of $\bk(x)$ defined by $\sigma(h):=g^*h$ for any $h\in \bk(x)$. 
\begin{que}[Picard-Vessiot problem]\label{quepvproblem}
	Does there exist a Picard-Vessiot extension of $\bk(x)$ for the
	linear difference equation 
	$$\sigma^l(y)-\sum_{i=0}^{l-1}a_{i}\sigma^i(y)=0.$$
	inside the ring of $\bk$-valued
	sequences?
\end{que}
For more background on difference equations and the aforementioned two problems, we refer the reader to \cite{Wibmer} and \cite{Put2003}.
It is shown in \cite[Theorem 1.5]{Ghioca2018} that Question \ref{quepvproblem} has positive answer when $\deg g\geq 2.$

\section{Known results and improvements}\label{knownresults}
The DML conjecture is only known in a few special cases.
Here are some notable results:
\begin{thm}\label{thmdmlknownex}The DML conjecture is true in the following cases:
\begin{points}
\item(Bell, Ghioca and Tucker, \cite{Bell2010}) The map $f$ is an \'etale endomorphism.
\item(Xie, \cite{Xie2017a}) The map $f$ is an endomorphism of $\A^2$ defined over $\overline{\Q}.$
\item(Ghioca and Xie, \cite{Ghioca2018}) The map $f$ is a skew-linear map on $\P^1\times \A^l$ taking form 
$$(x,y)\mapsto (g(x), A(x)y)$$ where $g$ is an endomorphism of $\P^1$ of degree $\geq 2$ and 
$A(x)$ is a matrix in $M_{l\times l}(k(x))$.
\end{points}
\end{thm}

The proof of (i) and (iii) relies on the $p$-adic interpolation lemma (\cite[Theorem 1]{Poonen2014} and \cite[Theorem 3.3]{Bell2010}).
Such $p$-adic interpolation method for studying DML backs to the proof of Skolem-Mahler-Lech Theorem \cite{Lech1953} and it was mainly developed by Bell, Ghioca and Tucker in a series of papers. In Section \ref{sectionautaff}, we reformulate and generalize this method in the language of Berkovich space.

The proof of (ii) is based on the theory of valuation tree introduced by Favre and Jonsson in \cite{Favre2004} and developed in \cite{Favre2007,Favre2011,Xie2015ring,Xie2017a}. It also relies on Siegel's integral point theorem \cite[Theorems 8.2.4 and 8.5.1]{Lang1983}
and the strategy developed in \cite{Xie2014}.

\medskip

More known results can be found in \cite{Bell2010,Benedetto2012,Ghioca2009,Ghioca2008,Xie2014,Ghioca2018,Ghioca2019a,Xie2017a,Bell2016}.
Moreover, a very good introduction of the DML conjecture can be found in the monograph \cite{Bell2016} of Bell, Ghioca and Tucker.

\medskip

Argument in \cite[Page 4]{Ghioca2018} shows that (iii) of Theorem \ref{thmdmlknownex} implies  Conjecture \ref{connonconsml} (=SML with non-constant coefficients) when $\deg g\geq 2.$

With a few modifications of the proof in \cite{Xie2017a}, one can generalize (ii) of Theorem \ref{thmdmlknownex} for any endomorphism of $\A_{\C}^2$.
\begin{thm}\label{thmendoplaneoverc}
The DML conjecture holds when $f$ is an endomorphism of $\A^2_{\C}.$
\end{thm}
Theorem \ref{thmendoplaneoverc} implies  Conjecture \ref{connonlinearsml} (= non-linear SML) when $l\leq 2$.
We will show how  to modify the proof in \cite{Xie2017a} to get Theorem \ref{thmendoplaneoverc} in Section \ref{sectiondmladc}.

\medskip

Using Theorem \ref{thmendoplaneoverc},  we may generate \cite[Theorem 0.3]{Xie2017a} to any endomorphisms over $\C$.

\begin{cor}\label{corextensionbyonepoly}
DML(1) holds for $f: \A^{N}\to \A^{N}$ over $\C$ taking form 
$$(x_1,\dots, x_{N})\mapsto (F_1(x_1), \dots, F_{N}(x_N)).$$
where $F_i$ are polynomials. 
\end{cor}

We can also get Corollary \ref{corextensionbyonepoly} by \cite[Theorem 0.3]{Xie2017a} and  a specialization argument as in \cite{Ghioca2017b}.

\subsection{DML for endomorphisms of $\A^2_{\C}$}\label{sectiondmladc}
\proof[Proof of Theorem \ref{thmendoplaneoverc}] In the proof of \cite[Theorem 1]{Xie2017a}(=(ii) of Theorem \ref{thmdmlknownex}), 
the assumption that $f$ is defined over $\overline{\Q}$ is only used for two purposes:

\begin{points}
\item[(1)]The first one is to apply Siegel's integral points theorem \cite{Hindry2000}, which asserts that for a curve $C$ defined over a number field $K$ whose normalization has at least $3$ boundary points, for every finite set of places $S$ of $\sM_K$ containing of archimedean places, the set $S$-integral points (for any models of $C$ over $S$)  is finite. 
\item[(2)] The second one is the Northcott property for Weil heights of $\overline{\Q}$ points. 
\end{points}
In the over $\C$ case, we can choose a subring $R$ of $\C$ which is finitely generated over $\Z$ such that the setups in Theorem \ref{thmendoplaneoverc} are defined over $R.$
Replace the above (1) and (2) by their versions for finitely generated fields, we may conclude the proof.

For (2), we replace Weil's height machine over number field to 
the height machine for the Moriwaki height over finitely generated fields \cite{Moriwaki2000}. The Moriwaki height satisfies all the properties we need in the proof, in particular, it has the Northcott property.

For (1), we replace Siegel's integral points theorem over number field to its version for $R$-points \cite[Theorem 4]{Lang1960}. 
The statement in \cite[Theorem 4]{Lang1960} is not exactly the one we need. However, it is not hard to deduce to the version we need as follows:
In \cite[Theorem 4]{Lang1960}, the curve is assumed to be smooth and we need to treat singular curves. This assumption is easy to get by taking normalization.
In \cite[Theorem 4]{Lang1960}, the curve is asked to have genus at least $1$. In our case, the curve $C$ may be $\P^1$ with at least three removed points.
In this case, we may take a suitable \'etale cover $C'$ of $C$ satisfying the assumption in \cite[Theorem 4]{Lang1960}. 
By the version of Hermite-Minkowski theorem over finitely generated ring
 \cite[Theorem 2.9]{Harada2009}, Siegel's integral points theorem for $C'$ implies Siegel's integral points theorem for $C.$
\endproof

\section{DML in positive characteristic}\label{dmlpositivec}
As shown in the following example, the DML conjecture does not hold in positive characteristic. 

\begin{exe}\cite[Example 3.4.5.1]{Bell2016}\label{exenotdml} Let $p$ be a prime and $\bk=\overline{\F_p(t)}$. Let $f: \A^2\to \A^2$ be the endomorphism defined by $(x,y)\mapsto (tx, (1-t)y).$
	Set $V:=\{x+y=1\}$ and $e=(1,1).$ Then $\{n\geq 0|\,\, f^n(e)\in V\}=\{p^n|\,\, n\geq 0\}.$
\end{exe}
More examples can be found in \cite{Ghioca2019a,Corvaja2021}.
So we have the following two natural questions when $\Char\, \bk>0$:
\begin{points}
\item[(a)] What is the form of the set
$\{n\geq 0|\,\, f^n(x)\in V\}$?
\item[(b)] DML holds for which $(X,f)$?
\end{points}

As an answer of (a), Ghioca  and Scanlon proposed a variant of the Dynamical Mordell-Lang conjecture in positive characteristic (=$p$-DML) in \cite[Conjecture 13.2.0.1]{Bell2016}, which asserts that up to a finite set, $\{n\geq 0|\,\, f^n(x)\in V\}$ is a finite union of arithmetic progressions along with finitely many 
sets taking form $$\{\sum_{i=1}^mc_ip^{k_in_i}|\,\, n_i\in \Z_{\geq 0}, i=1,\dots,m\},$$
where $m\in \Z_{\geq 1}, k_i\in \Z_{\geq 0}, c_i\in \Q$.

Surprisingly, the $p$-DML is even open for endomorphisms of semi-abelian varieties.
Indeed in \cite{Corvaja2021}, Corvaja, Ghioca, Scanlon and Zannier showed that $p$-DML for endomorphisms of semi-abelian varieties
is equivalent with some difficult diophantine problem on polynomial-exponential equations in characteristic $0$.
Partial results can be found in \cite{Ghioca2019a,Corvaja2021}. These works essential rely on Hrushovski's resolution of the Mordell-Lang conjecture in positive characteristic \cite{Hrushovski1996} and the further description of the intersection $\Gamma\cap V$ in \cite{Moosa2002,Moosa2004,Ghioca2008b, Ghioca2023}.

\medskip

For (b), the author in \cite{Xie2023} (and also Ghioca and Scanlon) guessed that DML holds for 
most dynamical systems and the counter-examples often involve
some group actions.
Outside the case of endomorphisms of semi-abelian varieties which was studied in \cite{Ghioca2019a,Corvaja2021}, only a few results are known:
\begin{thm}\label{thmdmlpochar}DML holds in the following cases:
\begin{points}
\item(Xie, \cite[Theorem A]{Xie2014}) The map $f: \A^2\to \A^2$ is a birational endomorphism with $\la_1(f)>1$.
\item(Xie, \cite[Theorem 1.4]{Xie2023}) The map $f: X\to X$ is an automorphism of a projective surface $X$ with $\la_1(f)>1$.
\item(Yang, \cite[Theorem 1.0.2]{Yang2023}) Let $\bk$ be a complete non-archimedian valuation field with $\Char\, \bk=p>0.$
Let $\bk^{\circ}$ be its valuation ring and $\bk^{\circ\circ}$ be the maximal ideal.
Assume that $f:\P^N_{\bk}\to \P^N_{\bk}$  is totally inseparable  and it is a lift of Frobenius i.e. it
takes the following form:
$$f: [x_0:\dots :x_N]\mapsto [x_0^q+g_0(x_0^p,\dots,x_N^p):\dots :x_N^q+g_N(x_0^p,\dots,x_N^p)]$$
where $q$ is a power of $p$ and $g_i$ are homogenous polynomials of degree $q/p$ in $\bk^{\circ\circ}[x_0,\dots,x_N].$
\end{points}
\end{thm}
In Theorem \ref{thmdmlpochar}, $\la_1(f)$ stands for the first dynamical degree of $f$ \cite{Russakovskii1997, Dinh2005, Truong2020,Dang2020,boucksomfavrejonsson} (c.f. Section \ref{sectionksmuupb}).
The proofs of (i) and (ii) indeed work in any characteristic. The situation of (iii) only appears in positive characteristic case. 

\subsubsection*{Idea to prove (i):} For the simplicity, assume that $\bk=\overline{\Q}.$ Let $x\in \A^2(\bk)$ be a non-periodic point and $V$ be an irreducible curve in $\A^2.$ First we find a ``good" compactification $X$ of $\A^2$ w.r.t $f$ in the sense of \cite{Favre2011}. We view $f$ as a birational self-map on $X.$
After replacing $f$ by a suitable iterate, there is a superattracting fixed point $q\in \partial X:=X\setminus \A^2$ such that $f$ maps $\partial X\setminus I(f)$ to $q$. By some geometric arguments, we reduce to the case where $q\in \overline{V}.$ Next, we prove a local version of DML near $q$ for the topology on $X(\bk)$ induced by any valuation on $k.$ In the end, we conclude the proof by the Norhcott property.

\subsubsection*{Idea to prove (ii):} For the simplicity, assume that $\bk=\overline{\Q}.$ Let $x\in X$ be a non-periodic point and $V$ be an irreducible curve in $X.$ Using an argument backs to \cite{Silverman1991,Kawaguchi2008}, we construct two canonical height functions $h^+$ and $h^-$
associated to two nef numerical classes $\theta^+,\theta^-\in N^1(X)_{\R}$ such that 
$$f^*h^+=\la_1(f)h^+ \text{ and } f^*h^-=\la_1(f)^{-1}h^-.$$
By Hodge index theorem, $\theta^++\theta^-$ is big and nef.
Easy to show that if $V\cdot (\theta^++\theta^-)=0$, then $V$ is $f$-periodic.
As an application of the Hodge index theorem, one can show that there exists $\sigma\in \Gal(\overline{\Q}/\Q)$ such that $\sigma(\theta^+)=\theta^-.$ So $V\cdot (\theta^++\theta^-)>0$ implies that $V\cdot \theta^->0.$
Considering $h^-(f^n(x))$ for those $f^n(x)\in V$, we conclude the proof by the Northcott property. 

\subsubsection*{Idea to prove (iii):}
For the simplicty, we assume that $\bk$ is a local field.
We extend $f$ to an endomorphism on the model $\P^N_{\bk^{\circ}}$.
Let $x$ be a point in $\P^N(\bk)$ and $V$ be an irreducible subvariety of dimension $\geq 1$.
Assume that $O_f(x)\cap V$ is Zariski dense in $V.$
Observe that in this case, $f$ is super-attracting at each $\bk$-point. Our idea is to mimic the process of constructing stable/unstable manifold in hyperbolic dynamics. More precisely, we show that at each closed point $\widetilde{x}$ in the special fiber, among all of its lifts
$x\in \P^N(\bk)$ there are some special ones which looks like the ``unstable" manifolds through $\widetilde{x}$. The unstable lifts of $\widetilde{x}$ form a Cantor set. To translate this idea from hyperoblic dynamics to our purely algebraic setting, one need to apply the techniques of the jet scheme and the critical scheme introduced in \cite{Roessler2013}. In this step, it is crucial to use the assumption that $f$ is totally inseparable.
Our assumption that $O_f(x)\cap V$ is Zariski dense in $V$ implies that for every closed point $\widetilde{x}$ in the special fiber of $V$, it has an unstable lift $x$ in $V(\bk)$. Pick a lift of $q$-Frobenius $\sigma$ in $\Gal(\overline{\bk}/\bk)$. The assumption that $f$ is a lift of Frobenius and the construction of the unstable lift $x$ shows that $\sigma(x)=f(x).$ It implies that $f(V)=\sigma(V).$ This implies that $V$ is periodic and concludes the proof. The idea of last step backs to Scanlon's proof of the dynamical Manin-Mumford conjecture for periodic points of lifts of Frobenius \cite{Xie2018}.


\section{Weak DML}\label{weakdml}
Let $X$ be a variety over a field $\bk$. Let $f: X\dashrightarrow X$ be a rational self-map.

\medskip

In \cite[Corollary 1.5]{Bell2015}, Bell, Ghioca and Tucker proved the Weak DML theorem, which assert that if the orbit of a point $x\in X(\bk)$ is Zariski dense, then
the  return time set $\{n\geq 0|\,\, f^n(x)\in V\}$ has zero density.  This version weakened the original DML in time.

\begin{thm}[Weak DML]\label{thmdml}
Let  $x$ be a point in $X_f(\bk)$ with $\overline{O_f(x)}=X.$ Let $V$ be a proper subvariety of $X$. Then $\{n\geq 0|\,\, f^n(x)\in V\}$ is of Banach density zero in $\Z_{\geq 0}$; i.e. for every sequence of interval $I_n, n\geq 0$ in $\Z_{\geq 0}$ with $\lim_{n\to \infty}\# I_n=+\infty$, we have 
$$\lim_{n\to \infty}\frac{\#(\{n\geq 0|\,\, f^n(x)\in V\}\cap I_n)}{\#I_n}=0.$$
\end{thm}
This result was proved in earlier works in \cite[Theorem 2.5.8]{Favre2000a}, \cite[Theorem D, Theorem E]{Gignac2014} in a different form for natural density. 
See  \cite[Theorem 2]{Petsche2015}, \cite[Theorem 1.10]{Bell2020} and \cite[Theorem 1.17]{Xie2023} for different proofs. 

\medskip

Here we state and prove a uniform version of Weak DML. 
\begin{thm}[Uniform Weak DML]\label{thmdmlweaku}
Let $X$ be a variety over a field $\bk$. Let $f: X\dashrightarrow X$ be a rational self-map. 
Let $V$ be a proper closed subset of $X$. 
For every $\epsilon>0$, there is $N\geq 1$ such that for every $x\in X_f(\bk)$ with Zariski dense orbit and every intervals $I$ in $\Z_{\geq 0}$ with $|I|\geq N,$
we have $$\frac{\#(\{n\geq 0|\,\, f^n(x)\in V\}\cap I)}{\#I}<\epsilon.$$
\end{thm}
\rem In this version, the sparsity of the return time set is uniform for the point $x.$
It seems that even the DML conjecture does not imply this Uniform Weak DML directly.
\endrem
As suggested in \cite{Xie2023}, a natural way to study the (Uniform) Weak DML problem is to use the ergodic theory for constructible topology. 
We also get a generalization of Theorem \ref{thmdml} (c.f. Theorem \ref{thmperiodicapp}) for partial orbits, in which we don't need the whole orbit of $x$ to be well defined.
The idea of using ergodic theory for Zariski topology backs to the works of Favre and Gignac \cite{Favre2000a, Gignac2014}. 
In \cite{Xie2023}, the author studies ergodic theory for constructible topology instead of Zariski topology.
This idea also played an important role in \cite{Favre2022} for studying the entropy of dynamics of Berkovich spaces.

\medskip

The proof of Theorem \ref{thmdmlweaku} is ineffective. It is natural to ask the following question:
\begin{que}\label{queexplicitwdml}
Is there an explicit way to compute  $N$ in Theorem \ref{thmdmlweaku} using $X,V$ and $\epsilon$?
\end{que}

\subsection{Constructible topology}
Let $X$ be a variety over a field $\bk$. 
Denote by $|X|$ the underling set of $X$ with the constructible topology; i.e. the topology on a  $X$ generated by the constructible subsets (see~\cite[Section~(1.9) and in particular (1.9.13)]{EGA-IV-I}).
In particular every constructible subset is open and closed.
This topology is finer than the Zariski topology on $X.$ Moreover $|X|$ is (Hausdorff) compact.

 Using the Zariski topology, on may define a partial ordering on $|X|$ by $x\geq y$ if and only if $y\in \overline{x}.$
The noetherianity of $X$ implies that this partial ordering satisfies the descending chain condition: for every chain in $|X|$, 
$$x_1\geq x_2\geq \dots$$
there is $N\geq 1$ such that $x_n=x_N$ for every $n\geq N.$
For every $x\in |X|$, the Zariski closure of $x$ in $X$ is 
$U_x:=\overline{\{x\}}=\{y\in |X||\,\, y\leq x\}$ which is open and closed in $|X|.$

\medskip

Denote by $\sM(|X|)$ the space of Radon measures on $X$ endowed with the weak-$\ast$ topology.
The following is \cite[Theorem 1.12]{Xie2023}.
\begin{thm}\label{thmRadon}Every $\mu\in \sM(|X|)$ takes form 
$$\mu=\sum_{i\geq 0}a_i\delta_{x_i}$$
where $\delta_{x_i}$ is the Dirac measure at $x_i\in X$, $a_i\geq 0$.
\end{thm}

Let $\sM^1(|X|)$ be the space of probability Radon measures on $|X|.$ 
Since $|X|$ is compact, $\sM^1(|X|)$ is compact.  It is also sequentially compact as showed in the following result. 
\begin{cor}\cite[Corollary 1.14]{Xie2023}.\label{corsmxsc}
The space $\sM^1(|X|)$ is sequentially compact.
\end{cor}

%

%
%
%
%

\subsection{Proof of the Uniform Weak DML}
Let $f: X\dashrightarrow X$ be a dominant  rational self-map. Set $|X_f|:=\cap_{n\geq 0}f|_{X\setminus I(f)}^{-n}(X\setminus I(f)).$
Because every Zariski closed subset of $X$ is open and closed in the constructible topology,  $|X_f|$ is a closed subset of $|X|.$
The restriction of $f$ to $|X_f|$ is continuous. We still denote by $f$ this restriction.

 \medskip

  Denote by $\sP(X,f)$ the set of $f$-periodic points in $|X_f|.$
 Theorem \ref{thmRadon} implies directly the following lemma.

 \begin{lem}\cite[Lemma 5.3]{Xie2023}\label{leminvm} If $\mu\in \sM^1(|X_f|)$ with $f_*\mu=\mu$, then there are $x_i\in \sP(X,f), i\geq 0$ and $a_i\geq 0, i\geq 0$ with $\sum_{i=0}a_i=1$
 such that  
 $$\mu=\sum_{i\geq 0}\frac{a_i}{\#O_f(x_i)}(\sum_{y\in O_f(x_i)}\delta_y)$$
 \end{lem}

\proof[Proof of Theorem \ref{thmdmlweaku}]If Theorem \ref{thmdmlweaku} is not correct, there is  $\epsilon>0$ and a sequence of points
 $x_m\in X_f(\bk)$ with Zariski dense orbit and  
 intervals $I_m:=\{a_m,\dots, b_m\}$ in $\Z_{\geq 0}$ with $|I_m|\to \infty,$
such that $$\frac{\#(\{n\geq 0|\,\, f^n(x_m)\in V\}\cap I_m)}{\#I_m}\geq \epsilon.$$
Set $$\mu_n:=\frac{1}{I_m}\sum_{n\in I_m}\delta_{f^n(x_m)}.$$

Let $\eta_X$ be the generic point of $X.$
Set $Y:=|X_f|\setminus (\cup_{y\in \sP(X,f)\setminus \{\eta_X\}}\cup_{n\geq 0}f^{-n}(U_y))$, which is closed in $|X_f|.$
We have  $\Supp\,\, \mu_n\subseteq  Y.$
After taking subsequence, we may assume that $\mu_m$ converges to a probability measure $\mu\in \sM^1(|X_f|)$ in the weak-$\ast$ topology.
Observe that $$f_*\mu_m-\mu_m=\frac{1}{I_m}(\delta_{f^{b_m+1}(x_m)}-\delta_{f^{a_m}(x_m)})\to 0$$
as $m\to \infty.$
So $f_*\mu=\mu.$ Since $\Supp\,\, \mu\subseteq  Y$ and $\eta_X=\sP(X,f)\cap Y$, by Lemma \ref{leminvm}, $\mu=\delta_{\eta_X}.$
Note that the character function $1_V$ of $V$ is continuous on $|X|.$ Then  we have 
$$\epsilon\leq \frac{\#(\{n\geq 0|\,\, f^n(x_m)\in V\}\cap I_m)}{\#I_m}=\int 1_{U_V}\mu_m\to \int 1_{U_V}\delta_{\eta_X}=0.$$
We get a contradiction.
\endproof

\subsection{Partial orbits}
Define $|X|^+$ to be the disjoint union of $|X|$ with a point $\phi$. 
We extend the partial ordering on $|X|$ to $|X|^+$ by setting $x\geq \phi$ for every $x\in |X|^+.$
For every $x\in |X|^+$ define $U_x^+:=\{y\leq x\}.$ When $x\in X, U_x^+=U_x\cup \{\phi\}.$
It is clear that $\sM(|X|^+)=\sM(X)\oplus \R_{>0}\delta_{\phi}.$
By Theorem \ref{thmRadon} and Corollary \ref{corsmxsc}, $\sM(|X|^+)^1$ is sequentially compact
and every Radon measure on $|X|^+$ is atomic. 

\medskip

Let $f^+: |X|^+\to |X|^+$ be the endomorphism such that
$f^+|_{I(f)}$ is the map $I(f)\to \{\phi\}$ and $f^+|_{X\setminus I(f)}=f|_{X\setminus I(f)}$.
Since $I(f)$ is open and closed in $|X|$, $f^+$ is continuous.
Moreover, for $f^+$ preserves the ordering. 

\medskip

Define $I(f)^+:=f^{-1}(\phi)$ and $I(f):=I(f)^+\setminus \{\phi\}.$ For every $n\geq 0$, define $X_f(n):=\cap_{i=0}^nf^{-i}(|X|)$. Then $X_f(n)$ is a constructible subset of $X$. 
It is clear that $I(g^+)=I(g).$ For every $x\in X_f(n), n\geq 0$, define $$\mu_n(x):=\frac{1}{n+1}\sum_{i=0}^n\delta_{f^i(x)}.$$
For for every $y\in \sP(X,f)$, let $r_y$ be the minimal period of $f$. Define $$\nu_y:=r_y^{-1}\sum_{i=0}^{r_y}\delta_{f^i(y)}.$$
Let $\sM_d(X,f)$ be the subset of $\sM(X)$ taking form $\sum_{0=1}^la_i\nu_{y_i}$ with $a_i>0$, $y_i,i=1,\dots, l$ is a strictly decreasing sequence of points in $\sP(X,f).$
Indeed, we have $l\leq \dim X$. Set $\sM_d^1(X,f):=\sM_d(X,f)\cap \sM^1(X,f).$

\medskip

As shown in the proof of \cite[Theorem 1.17]{Xie2023},
Theorem \ref{thmdml} has a measure theoretic reformulation, which assert that, for  any point in $x\in X_f(\bk)$ with $\overline{O_f(x)}=X$ and any sequence of intervals $I_n, n\geq 0$ in $\Z_{\geq 0}$ with $\lim_{n\to \infty}\# I_n=+\infty$, we have 
$$|I_n|^{-1}\sum_{i\in I_n} \delta_{f^i(x)}\to \delta_{\eta_X}$$
in weak-$\ast$ topology.

The following result generalize Theorem \ref{thmdml}, for which we don't need to take the whole orbit.

\begin{thm}\label{thmperiodicapp}
Let  $N_m\in \Z_{\geq 0}, m\geq 0$ be a sequence tending to $+\infty$ and $x_m\in X_f(N_m-1)$.
Assume that $$\mu_{x_m,N_m}\to \mu$$ as $m\to \infty$  in the weak-$\ast$ topology.
Then $\mu_{x_m,N_m}\in \sM_d^1(X,f)$ 
\end{thm}

\proof[Proof of Theorem \ref{thmperiodicapp}]
As $x_m\in X_f(N_m-1)$ for every $m\geq 0$, we have $\supp\, \mu\subseteq \sM^1(|X|).$
Observe that $$f^+_*\mu_{x_m,N_m}-\mu_{x_m,N_m}=\frac{1}{N_m+1}(\delta_{(f^+)^{N_m+1}(x_m)}-\delta_{x_m})\to 0$$
as $m\to \infty.$
So $f^+_*\mu=\mu.$

By Theorem \ref{thmRadon}, we may write $\mu=\sum_{i\geq 1}a_i\delta_{x_i}$
where $x_i\in X$ are distinct, $a_i>0$ and $\sum_{i\geq 1} a_i=1$. Since $\mu$ is $f^+$-invariant,
we may write $\mu=\sum_{i=1}^l b_i\nu_{y_i}$ where $l\in \Z_{\geq 1}\cup \{+\infty\}$,
$b_i>0$, $\sum_{i\geq 1}b_i=1$ and $y_i$ are distinct points in $\sP(X,f).$
 Define $Y_i:=\cup_{j=0}^{r_{y_i}-1}U_{f^j(y_i)}^+.$ We have $f^+(Y_i)\subseteq Y_i.$
 
 \medskip

We first treat the case where  $b_1,b_2>0$ and $Y_1\not\subseteq Y_2$ and $Y_2\not\subseteq Y_1.$ 
Set $U_1:=Y_1\setminus Y_2$ and $U_2:=Y_2\setminus Y_1.$ Note that $U_i, i=1,2$ are open and closed in $|X|.$
For $i=1,2$, we have $$\int 1_{U_i}\mu_m\to \int 1_{U_i}\mu\geq b_i.$$
So for $m\gg 0$, we have $\int 1_{U_i}\mu_m\geq b_i/2.$ So there is a minimal $s_m\in \{0,\dots, N_m-1\}$ such that 
$f^{s_m}(x_m)\in U_1\cup U_2.$ Since $U_1\cap U_2=\emptyset$, there is a unique $t_m\in \{1,2\}$ such that $f^{s_m}(x_m)\in U_{t_m}.$
We claim that $f^{i}(x_m)\not\in Y_2\cup Y_1$ for $i=0,\dots, s_m-1.$
Otherwise, $f^{i}(x_m)\subseteq (Y_2\cup Y_1)\setminus (U_1\cup U_2)=Y_1\cap Y_2.$ Hence $f^{s_m}(x_m)\subseteq f^{s_m-i}(f^i(x_m))\subseteq Y_1\cap Y_2$, which is a contradiction. So the claim holds. 
For every $i\geq s_m$, $f^i(x_m)=f^{i-s_m}(f^{s_m}(x_m))\in Y_{t_m}.$ 
So $f^i(x_m)\not\in U_{3-t_m}$ for $i\geq s_m.$ This implies that $$\int 1_{U_{t_m}}\mu_m=0,$$ which is a contradiction. 

Hence for every distinct  $i,j\in  \{1,\dots, l\}$  either
$Y_i\subseteq Y_j$ or $Y_j\subseteq Y_i$. By the Noetherianity, we may assume that  $$Y_l\subsetneq \dots\subsetneq Y_1.$$
After replacing each $y_i$ by a suitable element in its orbit, we may ask $y_1>\dots >y_l.$ This concludes the proof.
\endproof

\subsection{An example}
The limiting measure in Theorem \ref{thmperiodicapp} is much more complicate than the one for a whole orbit c.f. Theorem \ref{thmdml}.
In this section, we give an example to show that these more complicated measure are necessary. 

\medskip

Let $X:=(\P_{\C}^1)^l.$
Let $g_i: \P^1\to \P^1, i=1,\dots, l$  endomorphisms of degree $d_i\geq 2$.
Assume that $d_1<\dots<d_l$.
Let $f: X\to X$ be the map $f=g_1\times \dots\times g_l.$
For every $i=1,\dots,l$ consider a sequence of points 
$o_i(n), n\geq 0$ with $g_i(o_i(0))=o_i(0)$ and $g_i(o_i(n))=o_i(n-1)$ for $n\geq 1.$
Assume that $o_i(1)\neq o_i(0)$ for every $i.$
Let $b_0,\dots,b_l\geq 0$ with $\sum_{i=0}^l b_i=1.$ 
Let $y_i$ be the generic point of $\P^{l-i}\times\{o_{l-i+1}\}\times \dots\times \{o_l\}$ ($=(\P^1)^l$ when $i=0$).
We will construct a sequence of points $x_m, m\geq 1$ and $N_m\geq 1$ such that 
$\mu(x_m, N_m)\to \sum_{i=0}^lb_i\delta_{y_i}$ in weak-$\ast$ topology.

\medskip

Define $n_i(m):=\lceil m  \times b_i\rceil.$ 
Note that in $\P^{l}(\R),$ $[n_0(m):\dots: n_l(m)]$ tends to $[b_0:\dots: b_l]$ in the euclidien topology. 
Define $$x_m:=(o_0(\sum_{i=0}^{l-1}n_i(m)), o_1(\sum_{i=0}^{l-2}n_i(m)),\dots, o_{l-1}(n_0(m)))\in X=(\P_{\C}^1)^l$$
and $N_m:=\sum_{i=0}^ln_i(m).$
We will show that 
\begin{equation}\label{equmutenyi}\mu(x_m, N_m)\to \sum_{i=0}^lb_i\delta_{y_i}.
\end{equation}

By Corollary \ref{corsmxsc}, we only need to treat the case where $\mu(x_m, N_m)$ converges to a measure $\mu$ in $\sM^1(X).$
We do the proof by induction on $l.$
First assume that  $l=1$. 
Since $\int 1_{o_0}\mu(x_m, N_m)= n_0(m)/N_m\to b_0$, 
by Theorem \ref{thmperiodicapp}, $\mu(x_m, N_m)\to b_0\delta_{y_1}+b_1\delta_{y_0}$.
Now assume that $l\geq 2$ and (\ref{equmutenyi}) holds for $l-1$.
Define $Y:=(\P^1)^{l-1}$ and $h: Y\to Y$ by $g_2\times\dots\times g_l.$
Then $X=\P^1\times Y$ and $f=g_1\times h.$ Let $\pi_1, \pi_Y$ be the projection from $X$ to $\P^1$ and to $Y$ respectively.
Let $y_i', i=0,\dots,l$ be the generic point of $\P^{l-i-1}\times\{o_{l-i}\}\times \dots\times \{o_l\}$ in $Y$.
Let $z_0,z_1$ be the generic points of  $\P^1$ and $o_1(0)$ respectively.

The induction hypothesis implies that 
\begin{equation}\label{equpiomu}(\pi_1)_*\mu=(b_0+\dots+b_{l-1})\delta_{z_0} +b_l\delta_{z_1}
\end{equation}
and 
\begin{equation}\label{equpiymu} (\pi_Y)_*\mu=b_0\delta_{y_0'}+\dots +(b_{l-1}+b_l)\delta_{y_{l-1}'}.
\end{equation}
By \cite[Corollary 2.19 and Proposition 3.14]{Xie2022a}, $\mu$ takes form 
$$\mu=\sum_{i=0}^l c_i\delta_{z_{s(i)}\times y_{t(i)}'}$$
where $c_i\geq 0$, $\sum_{i=0}^l c_i=1$ and $s, t$ are increasing functions from $\{0,\dots,l\}$ to $\{0,1\}$ and to $\{0,\dots,l-1\}$ respectively satisfying 
$\max\{s(i)-s(i-1),t(i)-t(i-1)\}\geq 1$ for every $i=1,\dots,l.$
This follows that $(s(l), t(l))=(1,l-1)$ and $s^{-1}(1)\cap t^{-1}(l-1)=\{l\}.$
By (\ref{equpiomu}) and (\ref{equpiymu}), we have $c_l=b_l.$
By (\ref{equpiomu}), we get $s(i)=0$ for $i\in\{0,\dots, l-1\}$ with $b_i\neq 0.$
Then we concludes the proof by (\ref{equpiymu}).

\section{KSM's upper bound for arithmetic degree}\label{sectionksmuupb}
A notable application of  weak DML is \cite[Proposition 3.11]{Jia2021}, in which Jia, Shibata, Zhang and the author generalized  Kawaguchi-Silverman-Matsuzawa's upper bound=(KSM's upper bound) \cite[Theorem 1.4]{Matsuzawa2020a}
on arithmetic degree to singular case.  Recently, Song \cite[Theorem 3.2 and Theorem 4.6]{Song2023} gave a different and more conceptual proof of  the same upper bound using Arakelov geometry. Indeed Song's result is more general, which applies for higher dimensional arithmetic degree.
In Matsuzawa's original work \cite[Theorem 1.4]{Matsuzawa2020a}, he indeed proved a uniform version of the upper bound. 
In this section we use Theorem \ref{thmperiodicapp}(=Weak DML for partial orbit) to give a new proof of this the Uniform KSM's upper bound. Our result is indeed more general, it works on singular varieties for both number field and function field in any characteristic, moreover we don't need the whole orbit of $x$ is well-defined.
\begin{thm}[Uniform KSM's upper bound]\label{thmuniKSM'sineq}
Let $\bk$ be either $\overline{\Q}$ or the algebraically closure of a function field $K(B)$ for a smooth projective curve $B$ over an algebraically closed field $K.$
Let $f: X\dashrightarrow X$ be a dominant rational map defined over $\bk$. Let $h$ be any Weil height on $X$ associated to some ample line bundle.
We denote by $h^+:=\max\{h,1\}.$

Then for
any $\epsilon>0$, there exists $C>0$ such that
$$h^+(f^n(x))\leq C(\la_1(f)+\epsilon)^nh^+(x)$$
for all $n\geq 0$ and $x\in X_f(\bk, n)$. In particular, for any $x\in X_f(\bk)$, we have $$\overline{\alpha_f}(x)\leq \la_1(f),$$
where $\overline{\alpha_f}(x)$ stands for the upper arithmetic degree.
\end{thm}

Our proof here is based on the proof of the (non-uniform) KSM's upper bound by the author given in the course ``Topics in algberaic geometry: Arithmetic dynamics" in Peking university in 2022. 
The current uniform version of Theorem \ref{thmuniKSM'sineq} is inspirited by Matsuzawa's question, who asked whether the KSM's upper bound of arithmetic degree can be maked to uniform for singular varieties.

To prove Theorem \ref{thmuniKSM'sineq}, we first recall the definition and basic properties of dynamical degree and arithmetic degree.

\subsubsection{The dynamical degrees}\label{subsectiondydeg}
Let $X$ be a variety over a field $\bk$ and $f: X\dashrightarrow X$ a dominant rational self-map.
Let $X'$ be a normal projective variety which is birational to $X$.
Let $L$ be an ample (or just nef and big) divisor on $X'$.
Denote by $f'$ the rational self-map of $X'$ induced by $f$.

For $i=0,1,\dots,\dim X$, and $n\geq 0$,  let $(f'^n)^*(L^i)$ be the $(\dim X-i)$-cycle on $X'$ as follows: let $\Gamma$ be a normal projective variety with a birational morphism $\pi_1\colon\Gamma\to X'$ and a morphism $\pi_2\colon\Gamma\to X'$ such that $f'^n=\pi_2\circ\pi_1^{-1}$.
Then $(f'^n)^*(L^i):= (\pi_1)_*\pi_2^*(L^i)$.
The definition of $(f'^n)^*(L^i)$ does not depend on the choice of $\Gamma$, $\pi_1$ and $\pi_2$.
The $i$-th \textit{dynamical degree} of $f$ is
$$
\la_i(f):=\lim_{n\to\infty}((f'^n)^*(L^i)\cdot L^{\dim X-i})^{1/n}.
$$
The limit converges and does not depend on the choice of $X'$ and $L$
\cite{Russakovskii1997, Dinh2005, Truong2020,Dang2020}.
Moreover, if $\pi: X\dashrightarrow Y$ is a generically finite and dominant rational map between varieties and $g\colon Y\dashrightarrow Y$ is a rational self-map such that $g\circ\pi=\pi\circ f$, then $\la_i(f)=\la_i(g)$ for all $i$. This can be shown by combing \cite[Theorem 1]{Dang2020} with the projection formula.
Another way is to proof it is to apply the product formula for relative dynamical degrees (c.f. \cite{Dang} and \cite[Theorem 1.3]{Truong}) directly.

Recall the following basic property.
\begin{pro}\cite[Proposition 3.2]{Jia2021}\label{p:dyn_sub_var}
	Let $X$ be a variety over $\bk$ and $f\colon X\dashrightarrow X$ a dominant rational self-map.
	Let $Z$ be an irreducible subvariety in $X$ which is not contained in $I(f)$ such that $f|_Z$ induces a dominant rational self-map of $Z$.
	Then $\la_i(f|_Z)\leq \la_i(f)$ for $i=0,1,\dots,\dim Z$.
\end{pro}

\medskip

\subsubsection{Arithmetic degree}
The arithmetic degree was defined in \cite{Kawaguchi2016} over a number field or a function field of characteristic zero.
As said in \cite[Remark 1.14]{Matsuzawa2020a}, this definition can be extended to characteristic positive.
Here we follows the way \cite{Xie2023} the define the arithmetic degree in the general case.

Let $\bk=\overline{K(B)}$, where $K$ is an algebraically closed field and $B$ is a smooth projective curve.

Let $X$ be a normal and projective variety over $\bk.$ 
For every $L\in \Pic(X)$, we denote by $h_L: X(\bk)\to \R$ a Weil height associated to $L$ and the function field $K(B)$. It is unique up to adding a bounded function. 

%

\medskip

As in \cite{Jia2021,Xie2023}, we define an \textit{admissible triple} to be $(X,f,x)$ where $X$ is a quasi-projective variety over $\bk$, $f\colon X\dashrightarrow X$ is a dominant rational self-map and $x\in X_f(\bk)$.

We say that $(X,f,x)$ \textit{dominates} (resp.~\textit{generically finitely dominates}) $(Y,g,y)$ if there is a dominant rational map (resp.~generically finite and dominant rational map) $\pi\colon X\dashrightarrow Y$ such $\pi\circ f=g\circ\pi$, $\pi$ is well defined along $O_f(x)$ and $\pi(x)=y$.

We say that $(X,f,x)$ is \textit{birational} to $(Y,g,y)$ if there is a birational map $\pi\colon X\dashrightarrow Y$ such $\pi\circ f=g\circ\pi$ and if there is a Zariski dense open subset $V$ of $Y$ containing $O_g(y)$ such that $\pi|_U: U:=\pi^{-1}(V)\to V$ is a well-defined isomorphism and $\pi(x)=y$.
In particular, if $(X,f,x)$ is birational to $(Y,g,y)$, then $(X,f,x)$ generically finitely dominates $(Y,g,y)$.

\begin{rem}
	\leavevmode
	\begin{enumerate}
		\item If $(X,f,x)$ dominates $(Y,g,y)$ and if $O_f(x)$ is Zariski dense in $X$, then $O_g(y)$ is Zariski dense in $Y$.
		Moreover, if $(X,f,x)$ generically finitely dominates $(Y,g,y)$, then $O_f(x)$ is Zariski dense in $X$ if and only if $O_g(y)$ is Zariski dense in $Y$.
		\item Every admissible triple $(X,f,x)$ is birational to an admissible triple $(X',f',x')$ where $X'$ is projective.
		Indeed, we may pick $X'$ to be any projective compactification of $X$, $f'$ the self-map of $X'$ induced from $f$, and $x'=x$.
	\end{enumerate}
\end{rem}


As in \cite{Jia2021,Xie2023}, we will associate to an admissible triple $(X,f,x)$ a subset $$A_f(x)\subseteq [1,\infty].$$
Indeed we have $A_f(x)\subseteq [1,\la_1(f)]$ by Theorem
\ref{thmuniKSM'sineq}.

We first define it when $X$ is projective. Let $L$ be an ample divisor on $X$, we define
$$A_f(x)\subseteq [1,\infty]$$
to be the limit set of the sequence $(h_L^+(f^n(x)))^{1/n}$, $n\geq 0$, where $h_L^+(\cdot):=\max\{h_L(\cdot),1\}$.

The following lemma is \cite[Lemma 2.7]{Xie2023}, whose proof is the same as \cite[Lemma 3.8]{Jia2021}.
It shows that the set $A_f(x)$ does not depend on the choice of $L$ and is invariant in the birational equivalence class of $(X,f,x)$.

\begin{lemma}\cite[Lemma 2.7]{Xie2023}\label{lemsingwilldef}
	Let $\pi\colon X\dashrightarrow Y$ be a dominant rational map between projective varieties.
	Let $U$ be a Zariski dense open subset of $X$ such that $\pi|_U\colon U\to Y$ is well-defined.
	Let $L$ be an ample divisor on $X$ and $M$ an ample divisor on $Y$.
	Then there are constants $C\geq 1$ and $D>0$ such that for every $x\in U$, we have
	\begin{equation}\label{equationdomineq1}
		h_M(\pi(x))\leq Ch_L(x)+D.
	\end{equation}
	
	Moreover if $V:=\pi(U)$ is open in $Y$ and $\pi|_U\colon U\to V$ is an isomorphism, then 
	there are constants $C\geq 1$ and $D>0$ such that for every $x\in U$, we have
	\begin{equation}\label{equationbirdomineq}
		C^{-1}h_L(x)-D\leq h_M(\pi(x))\leq Ch_L(x)+D.
	\end{equation}
\end{lemma}

Now for every admissible triple $(X,f,x)$, we define $A_f(x)$ to be $A_{f'}(x')$ where $(X',f',x')$ is an admissible triple which is birational to $(X,f,x)$ such that $X'$ is projective.
By Lemma~\ref{lemsingwilldef}, this definition does not depend on the choice of $(X',f',x')$.

\medskip

\subsubsection{The arithmetic degree.}\label{subsec_arithdeg}
We define (see also \cite{Kawaguchi2016}):
\[
\overline{\alpha}_f(x):=\sup A_f(x),\qquad\underline{\alpha}_f(x):=\inf A_f(x).
\]
And call them \emph{upper/lower arithmetic degree}.
We say that $\alpha_f(x)$ is well-defined and call it the \textit{arithmetic degree} of $f$ at $x$, if $\overline{\alpha}_f(x)=\underline{\alpha}_f(x)$;
and, in this case, we set
\[
\alpha_f(x):=\overline{\alpha}_f(x)=\underline{\alpha}_f(x).
\]
By Lemma~\ref{lemsingwilldef}, if $(X,f,x)$ dominates $(Y,g,y)$, then $\overline{\alpha}_f(x)\geq \overline{\alpha}_g(y)$ and $\underline{\alpha}_f(x)\geq\underline{\alpha}_g(y)$.

Applying Inequality~\eqref{equationdomineq1} of Lemma~\ref{lemsingwilldef} to the case where $Y=X$ and $M=L$, we get the following trivial upper bound:
let $f\colon X\dashrightarrow X$ be a dominant rational self-map, $L$ any ample line bundle on $X$ and $h_L$ a Weil height function associated to $L$;
then there is a constant $C\geq 1$ such that for every $x\in X\setminus I(f)$, we have
\begin{equation}\label{equationtrivialupper}
	h_L^+(f(x))\leq Ch_L^+(x).
\end{equation}
For a subset $A\subseteq [1,\infty)$, define $A^{1/\ell}:= \{a^{1/\ell}\mid a\in A\}$.

We have the following simple properties, where the second half of \ref{eq:alpha_pow} used Inequality~\eqref{equationtrivialupper}.
\begin{pro}\label{probasicaf}We have:
	\begin{enumerate}
		\item $A_f(x)\subseteq [1,\infty)$.
		\item $A_f(x)=A_f(f^{\ell}(x))$, for any $\ell\geq 0$.
		\item \label{eq:alpha_pow}
		$A_{f}(x)=\bigcup_{i=0}^{\ell-1}(A_{f^{\ell}}(f^i(x)))^{1/\ell}$.
		In particular, $\overline{\alpha}_{f^{\ell}}(x)=\overline{\alpha}_{f}(x)^{\ell}$, $\underline{\alpha}_{f^{\ell}}(x)=\underline{\alpha}_{f}(x)^{\ell}$.
	\end{enumerate}
\end{pro}

\proof[Proof of Theorem \ref{thmuniKSM'sineq}(=the Uniform KSM's upper bound)]
Fix an ample line bundle $A$ on $X$. Let $h_A$ be a Weil height associated with $A.$
Assume that $h_A\geq 1$. In this case, we have $h^+=h_A.$
By Lemma \ref{lemsingwilldef}, after taking normalization, we may assume that $X$ is normal.

\medskip
We first need two height inequalities.
The first one is by Lemma \ref{lemsingwilldef}.
There is $C_1>0$ such that for every $X(\bk)\setminus I(f)$, we have
\begin{equation}\label{equheighttrivial}
	h_A(f(x))\leq C_1h_A(x).
\end{equation}
This is a weak estimation which holds on a large domain $X(\bk)\setminus I(f).$

\medskip

We next prove a stronger inequality, which holds only on some Zariski open subset of $X(\bk)\setminus I(f).$
\begin{lem}\label{lemstronginsiu}For every $\psi>\la_1(f)$, there is $N=N(\psi)\geq 1$ and $C_2=C_2(\psi)>0$ and a proper closed subset $B=B(\psi)$ of $X$ such that 
for every $n\geq 1$ and 
$x\in X_f(Nn)(\bk)$ satisfying $f^{Ni}(x)\not\in B$ for $i=0,\dots, n-1$, then we have 
$$h_A(f^{Nn}(x))\leq C_2\psi^{Nn}h(x).$$
\end{lem}
The proof of Lemma \ref{lemstronginsiu} is a simple application of Siu's inequality, we postpone the proof to the end of this section.

Combining Lemma \ref{lemstronginsiu} and (\ref{equheighttrivial}), we get that for every 
$\psi>\la_1(f)$, there is $C_3=C_3(\psi)>0$ and a proper closed subset $B_1=B_1(\psi)$ of $X$ such that 
for every $n\geq 1$ and 
$x\in X_f(n)(\bk)$ satisfying $f^{i}(x)\not\in B_1$ for $i=0,\dots, n-1$, then we have 
\begin{equation}\label{equgoodbounsiu}h_A(f^{n}(x))\leq C_3\psi^{n}h(x).
\end{equation}

\medskip

We do the proof of Theorem \ref{thmuniKSM'sineq} by induction on the dimension of $X.$ 
When $\dim X=0$, nonthing to prove. Now assume that $\dim X\geq 1.$
By (\ref{equheighttrivial}), if Theorem \ref{thmuniKSM'sineq} holds for a positive iterate of $f$, then it holds for $f.$
Assume that Theorem \ref{thmuniKSM'sineq} does not hold for our $(X,f).$

Then there are  $\epsilon>0$,  a sequence $N_m\geq 1$ satisfying $N_m\to \infty$, a sequence of points
$x_m\in X_f(N_m-1)(\bk)$ such that 
\begin{equation}\label{equcontra}h_A(f^{N_m}(x_{m}))\geq (\la_1(f)+\epsilon)^{N_m}h_A(x_m).
\end{equation}
After taking subsequence, we may assume that $\mu(x_m,N_m)\to \mu$ for some $\mu\in \sM^1(X).$
By theorem \ref{thmperiodicapp}, $\mu\in \sM_d^1(X)$.
We may write $$\mu=\sum_{i=0}^la_i\nu_{y_i}$$ with $a_0\geq 0,$ $a_i>0$ for $i\geq 1$, $y_1=\eta_X$ and  $y_i,i=1,\dots, l$ is a strictly decreasing sequence of points in $\sP(X,f)$ with $l\leq \dim X.$
Set $U_{1}:=\cup_{i=0}^{r_{y_1}-1}U_{y_1}$ if $l\geq 1$, otherwise set $U_1=\emptyset.$
Set $V:=X\setminus U_{1}$.

For every $m\geq 0$, there is a largest $s_m\geq 0$ such that 
$f^s_m(x_m)\in V.$ Hence $f^i(x_m)\in V$ for $i=0,\dots, s_m$ and $f^i(x_m)\not\in V$ for $i> s_m$.
Note that $U_1$ is an $f$-invariant Zariski closed subset of $\dim U_1<\dim X.$
Pick $\psi\in (\la_1(f),\la_1(f)+\epsilon).$
By Proposition \ref{p:dyn_sub_var}, the dynamical degree of restriction of $f^{r_{y_1}}$ to every irreducible component of $U_1$
is $\leq \la_1(f)^{r_{y_1}}.$
By the induction hypothesis for $f^{r_{y_1}}$ for every irreducible component of $U_1$  and (\ref{equheighttrivial}),
there is $C_3>0$ such that 
\begin{equation}\label{equinduarthiu}h_A(f^{N_m}(x_m))\leq C_3\psi^{N_m-s_m-1}h_A(f^{s_m+1}(x_m)).
\end{equation}

Since $$\int 1_V\mu(x_m,N_m)\to \int 1_V\mu=a_0,$$ we have 
$$(s_m+1)/(N_m+1)\to a_0$$ as $m\to \infty.$

\medskip

If $a_0=0$, by (\ref{equcontra}), (\ref{equheighttrivial}) and (\ref{equinduarthiu}) we get 
$$(\la_1(f)+\epsilon)^{N_m}h_A(x_m)\leq h_A(f^{N_m}(x_m))\leq C_3\psi^{N_m-s_m-1}C_1^{s_m+1}h_A(x_m).$$
Then we have
$$\log(\la_1(f)+\epsilon)\leq \frac{N_m-s_m-1}{N_m}\log \psi+\frac{1}{N_m}\log C_3+\frac{s_m+1}{N_m}\log C_1.$$
Taking $m\to \infty$, we get a contradiction. 

Now assume that $a_m>0$, then $s_m\to \infty.$
Set $W_m:=\{i=0,\dots, s_m| f^i(x_m)\in B_1\}$ and $w_m:=\#W_m.$
We have 
\begin{equation}\label{equwmzerod}\frac{w_m}{s_m+1}=\frac{w_m}{N_m+1}\times \frac{N_m+1}{s_m+1}=\frac{\int1_{B_1}\mu_m}{\int1_{V}\mu}\to 0/a_0=0.
\end{equation}
By (\ref{equheighttrivial}) and (\ref{equgoodbounsiu}), the argument in \cite[Lemma 3.12]{Jia2021} implies that 
$$h_A(f^{s_m}(x_m))\leq C_3^{w_m+1}C_1^{2w_m}\psi^{s_m-w_m}h_A(x_m).$$
Then by (\ref{equinduarthiu}) and (\ref{equcontra}),
we get
$$(\la_1(f)+\epsilon)^{N_m}h_A(x_m)\leq h_A(f^{N_m}(x_m))\leq C_3\psi^{N_m-s_m-1}C_3^{w_m+1}C_1^{2w_m}\psi^{s_m-w_m}h_A(x_m).$$
Let $m\to \infty,$ we get a contradiction by (\ref{equwmzerod}).
\endproof

\proof[Proof of Lemma \ref{lemstronginsiu}]
For every $n\geq 1$, denote by $X_n$ the normalization of the graph of $f^n: X\dashrightarrow X.$
Denote by $\pi_n: X_n\to X$ the first projection and $f_n: X_n\to X$ the second projection.
We have $f_n=f^n\circ \pi_n.$
By Siu's inequality \cite[Theorem 2.2.13]{Lazarsfeld}, we have
$$f_n^*A\leq d_X\frac{\deg_A(f^n)}{A^{d_X}}\pi_n^*A$$
where $d_X=\dim X$ and $\leq$ means that the different is pseudo-effective.
Set $Q:=(d_X+1)/A^{d_X}$. Then there is a positive $\R$-divisor $B_n$ such that 
$$Q\deg_A(f^n)\pi_n^*A=f_n^*A+B_n.$$
There is $N\geq 0$ such that $Q\deg_A(f^N)\leq \psi^N.$
Hence for every $x\in X_f(N)(\bk)\setminus B_N$, we have 
$$\psi^Nh_A(x)+D\geq h_A(f^N(x))$$
 for some $D\geq 0.$
Then we have $$\psi^N(h_A(x)+D/(\psi^N-1))\geq h_A(f^N(x))+D/(\psi^N-1).$$
Then for every $n\geq 1$ and 
$x\in X_f(Nn)(\bk)$ satisfying $f^{Ni}(x)\not\in B$ for $i=0,\dots, n-1$, then we have 
$$h_A(f^{Nn}(x))\leq h_A(f^{Nn}(x))+D/(\psi^N-1)$$
$$\leq \psi^{Nn}(h_A(x)+D/(\psi^N-1))\leq (D/(\psi^N-1)+1)\psi^{Nn}h_A(x).$$
This concludes the proof.
\endproof

\section{Almost DML}\label{almostdml}
Another way to weaken the DML is in space, which means that we only ask it for almost all points in $X(\bk).$
In \cite{Xie2022}, the author proved a result in this direction. In this result, ``almost all points" is made to be precise using the adelic topology.

Assume that $\bk$ is an algebraic closed field of characteristic $0$ with finite transcendence degree over $\Q$.
In \cite[Section 3]{Xie2022}, the author has introduced the adelic topology on $X(\bk).$
It has the following basic properties (cf.~\cite[Proposition~3.18]{Xie2022}):
\begin{points}
		\item It is stronger than the Zariski topology.
		\item It is ${\mathsf{T}}_1$,
		i.e., for every distinct points $x, y \in X(\bk)$
	there are adelic open subsets $U, V$ of $X(\bk)$ such that
	$x \in U, y \notin U$ and $y\in V, x \notin V$.
	\item Morphisms between algebraic varieties over $\bk$ are continuous for the adelic topology.
	\item Flat morphisms are open with respect to the adelic topology.
	\item \label{ppty-of-adelic-topo:irreducible}
	The irreducible components of $X(\bk)$ in the Zariski topology
	are the irreducible components of $X(\bk)$ in the adelic topology.
	\item Let $K$ be any subfield of $\bk$ which is finitely generated over $\mathbb{Q}$
	and such that $X$ is defined over $K$ and $\overline{K} = \bk$.
	Then the action
	\[
	\Gal(\bk/K)\times X(\bk)\to X(\bk)
	\]
	is continuous with respect to the adelic topology.
\end{points}
	
	\begin{rem}\label{rem:intersection_of_finite_adelic_open}
		When $X$ is irreducible,
		the property (\ref{ppty-of-adelic-topo:irreducible}) above implies that
		the intersection of finitely many nonempty adelic open subsets of $X(\bk)$ is nonempty.
		So, if $\dim X\geq 1$, the adelic topology is not Hausdorff.
		In general, the adelic topology is strictly stronger than the Zariski topology.
	\end{rem}
	
	An impotent example of adelic open subsets is as follows:
	Let $L$ be a subfield of $\bk$
	such that its algebraic closure ${\overline{L}}$ is equal to $\bk$,
	$L$ is finitely generated over $\Q$,
	and $X$ is defined over $L$,
	i.e., $X=X_L\otimes_L \bk$ for some variety $X_L$ over $L$.
	Fix any embedding $\tau\colon L\hookrightarrow \C_p$ (resp. $\C$).
	Then, given any open subset $U$ of $X_L(\C_p)$ for the $p$-adic (resp.~Euclidean) topology,
	the union $X_L(\tau, U):= \cup_\iota \Phi_\iota^{-1}(U)$
	for all embeddings $\iota\colon \bk \to \C_p$
	extending $\tau$ is,
	by definition, an open subset of $X(\bk)$ in the adelic topology.
	Moreover $X_L(\tau, U)$ is empty if and only if $U=\emptyset$.

\medskip

As defined in \cite[Section 1.2.1]{Xie2022}, we say that a property holds for an adelic general point, if it holds on a non-empty adelic open subset.
The following result is \cite[Proposition 3.27]{Xie2022}.
\begin{pro}[Almost DML]\label{proadelicgedml}
DML holds for an adelic general point in $X(\bk)$.
\end{pro}
This proof of Almost DML is based on the $p$-adic interpolation lemma (see Theorem \ref{thminteraffinoid} for a slightly generalization and reformulation in the language of Berkovich space). This strategy backs to \cite{Amerik} and \cite{Bell2010}.
This Almost DML is find to be very useful in the recent studies of Zariski dense orbit conjecture and Kawaguchi Silverman conjecture, see  \cite{Xie2022,Jia2021,Matsuzawa}.

\section{Automorphisms on $\bk$-affinoid spaces}\label{sectionautaff}
The $p$-adic method, in particular the $p$-adic interpolation Lemma plays important role in the study of the DML conjecture \cite{Bell2016, Bell2010}. It also useful for studying other problems such as the Zariski dense orbit conjecture \cite{Amerik2008, E.Amerik2011, Xie2019} and the group of birational self-maps \cite{Cantata,Cantat2020}.
In this section, we explain this method in the language of Berkovich spaces \cite{Berkovich1990}.
This method indeed works for any complete non-archimedean field of characteristic $0$. However, as I known, all of its applications in the  DML conjecture is to apply it over some $p$-adic field. For this reason, we continue to use the name ``$p$-adic method". At the moment, there is no analogy of this method in positive characteristic. This is one of the difficulty for DML in positive characteristic.

\medskip

Denote by $\bk$ a complete valued field with a non-archimedean norm $|\cdot |.$ 
Denote by $\bk^{\circ}:=\{f\in \bk|\,\,|f|\leq 1\}$ the valuation ring and $\bk^{\circ\circ}:=\{f\in \bk|\,\,|f|< 1\}$ its maximal ideal. Denote by $\widetilde{\bk}:=\bk^{\circ}/\bk^{\circ\circ}$ the residue field.
Let $A$ be a $\bk$-affinoid algebra. Let $\|\cdot\|$ be a norm on $A$ and let $\rho(\cdot)$ be the
spectral seminorm on $A$.  We have $\|g\|\geq \rho(g)$ for all $g\in A.$ If $A$ is reduced, these two norms are equivalent to each other.

Let $f:\sM(A)\to \sM(A)$ be the endomorphism induced by an endomorphism $f^*:A\to A.$
Then $f$ induces an action of $\Z_{\geq 0}$ on $\sM(A)$.

The difference operator $\Delta_f:=f^*-\id$ is a  bounded linear operator on the Banach $\bk$-space. Write $\|\Delta_f\|$ the operator norm with respect to $\|\cdot\|$. Denote by $\rho(\Delta_f)$ the spectral of the operator $\Delta$.
We note that for any norm $\|\cdot\|$ on $A$, we have $\|\Delta_f\|\geq \rho(\Delta_f)$. 
The following lemma shows that when  $\rho(\Delta_f)<1$,  $f$ is an automorphism.
\begin{lem}\label{lemdflesoauto}If $\rho(\Delta_f)<1$ ( resp. $\|\Delta_f\|<1$), then $f$ is an automorphism. Moreover, we have 
$\rho(\Delta_f)=\rho(\Delta_{f^{-1}})$ ( resp. $\|\Delta_f\|=\|\Delta_{f^{-1}}\|$).
\end{lem}
\proof[Proof of Lemma \ref{lemdflesoauto}]We only prove it for $\rho(\cdot)$. The proof for $\|\cdot\|$ is similar. 
Denote by $g^*:A\to A$ the operator defined by $$g^*:=\sum_{i=0}^{\infty}(-1)^i\Delta_f^i.$$ Since $\rho(\Delta_f)<1$, the above series converges and $\rho(g^*-\id)=\rho(\Delta_f).$
We may check that $g^*\circ f^*=f^*\circ g^*=\id.$ Then $f$ is an automorphism and $\rho(\Delta_{f^{-1}})=\rho(g^*-\id)=\rho(\Delta_f).$ 
\endproof

\medskip



\subsection{Interpolation of iterates}
In this section, assume that  $\Char\, \bk=0.$
Define $R(\bk):=1$ if $\Char \,\widetilde{k}=0$ and $R(\bk):=p^{-\frac{1}{p-1}}$ if $\Char \,\widetilde{\bk}=p>0.$ 
This constant was introduced by Poonen in \cite{Poonen2014}, to study the interpolation of iterates of analytic self-maps of the $p$-adic polydiscs.

\medskip

We note that $R(k)\in (0,1]$ and $R(k)=1$ if and only if $\Char \,\widetilde{k} =0.$
It is easy to see that $|i!|\geq R(k)^i$ for all $i\geq 0.$

Denote by $\D:=\sM(\bk\{T\})$ the unit disc with the group structure given by the addition. 
There exists a natural embedding $\Z_{\geq 0}\subseteq \bk^{\circ}\subseteq \D(\bk)\subseteq \D.$
The following theorem generalized \cite[Lemma 4.2]{Bell2006}, \cite[Theorem 3.3]{Bell2010} and \cite[Theorem 1]{Poonen2014}.
Our proof is basically the same as the proof in \cite[Theorem 1]{Poonen2014}. 

\begin{thm}[=$p$-adic interpolation Lemma]\label{thminteraffinoid}If $\rho(\Delta_f)<R(\bk)$, then there exists a unique action $(\D,+)$ on $\sM(A)$ which extends the action of $\Z_{\geq 0}$ on $\sM(A)$.
\end{thm}

\rem
We note that if $\|\Delta_f\|<R(\bk)$ for any norm $\|\cdot\|$ on $A$, then we have $\rho(\Delta_f)<R(\bk)$. 
\endrem

\proof[Proof of Theorem \ref{thminteraffinoid}]
The uniqueness comes from the fact that $\Z_{\geq 0}\subseteq \D$ is Zariski dense.
We only need to prove the part of existence.

For any $h\in A$, we denote by $G(T,h)$ the analytic function in $\bk\{T\}\hotimes A$
by $$G(T,h):=\sum_{i\geq 0}\binom{T}{i}\Delta_f^i(h)=\sum_{i\geq 0}\frac{T(T-1)\dots(T-i+1)}{i!}\Delta_f^i(h).$$
It always converges.
Indeed, we have $|i!|\geq  R(\bk)^i$ and $\|\Delta_f^i(h)\|\leq \|\Delta_f^i\|\|h\|$.  Pick $R$ such that 
$\rho(\Delta_f)<R<R(\bk)$. For $i$ large enough, we have $\|\Delta_f^i\|<R^i.$
For $i$ large enough, we have $\|\frac{T(T-1)\dots(T-i+1)}{i!}\Delta_f^i(h)\|\leq (R/R(\bk))^i$. Then $G(T,h)$ converges. 

We define a linear map $\Phi^*:A\to \bk\{T\}\hotimes A$ by $h\mapsto G(T,h)$.
The above argument shows that $\|\Phi\|\leq \sup_{i\geq 1}\frac{\|\Delta_f^i\|}{i!}.$ Then $\Phi$ is bounded. 

We note that for any $n\in \Z_{\geq 0}$, we have 
\begin{equation}\label{equaboundgnh}
G(n,h)=\sum_{i\geq 0}\binom{n}{i}\Delta_f^i(h)=(\id+\Delta_f)^n(h)=(f^*)^n(h).
\end{equation}
Since $\Z_{\geq 0}\subseteq \D$ is Zariski dense in $\D$ and for any $n\in \Z_{\geq 0}$, we have $$G(n,1)=1, \text{ and } G(n,h_1h_2)=G(n,h_1)G(n,h_2)$$ for $h_1,h_2\in A,$
we get $$G(T,1)=1, \text{ and } G(T,h_1h_2)=G(T,h_1)G(T,h_2)$$ for $h_1,h_2\in A.$
It implies that $\Phi^*$ is a morphism of $\bk$-affinoid algebra.
It defines a morphism $\Phi:\D\times \sM(A)\to \sM(A).$

We only need to show that $\Phi$ is an action of a group. 
Since $\Z_{\geq 0}\subseteq \D$ is Zariski dense, by  (\ref{equaboundgnh}), $\Phi$ is an action of a semigroup.  Since $$\Phi(0,\cdot)=\id$$ and $\D$ is a group, $\Phi$ is a group action.
\endproof

\begin{cor}\label{corfiniteorder}If $\rho(\Delta_f)<R(\bk)$ and $f^l=\id$ for some $l\geq 1$, then $f=\id.$ 
\end{cor}
\proof[Proof of Corollary \ref{corfiniteorder}]Theorem \ref{thminteraffinoid} shows that the action of $f$ extends to an action of $\D$ on $\sM(A).$
For any point $t\in \D(\bk)$, we denote by $g_t$ the action of $t$. Since $f^l=\id$, we have $$g_{nl}=f^{nl}=\id$$ for all $n\geq 0.$
Since $l\Z_{\geq 0}$ is Zariski dense in $\D$, $g_t=\id$ for all $t\in \D(\bk).$ In particular, we get $f=g_1=\id.$
\endproof

\begin{rem}\label{remvecautosmall}
Assume $\rho(\Delta_f)<R(\bk)$ and denote by $\Phi:\D\times \sM(A)\to \sM(A)$ the action of $\D$ defined in Theorem \ref{thminteraffinoid}.

The natural morphism, $\bk\{T\}\to \bk[\epsilon]/(\epsilon^2)$ sending $T\to \epsilon$ induces a morphism 
$$\Phi^*\hotimes \bk[T]/(T^2):A\to \bk[\epsilon]/(\epsilon^2)\hotimes A=A\oplus \epsilon A.$$
It can be written explicitly by 
$$h\mapsto h+\epsilon \theta_f(h)$$
where $\theta_f:=\log(\id+\Delta_f)=\sum_{i\geq 1}\frac{(-1)^{i-1}}{i}\Delta_{f}^i.$
Since $\Phi^*\hotimes \bk[T]/(T^2)$ is a morphism of $\bk$-affinoid algebra,  $\theta_f$ is a vector field i.e. it satisfy the Leibniz's rule.
Such a construction was used in \cite{Cantata} for the proof of a birational version of Zimmer conjecture.
\end{rem}

The following results shows that even when Theorem \ref{thminteraffinoid} does not apply for $f$ directly, it may apply for a suitable iterate of $f.$
\begin{pro}\label{prodeltaiteratesmall}If $\|\Delta_f\|<1$ (resp. $\rho(\Delta_f)< 1$), then there exists $N$ such that $\|\Delta_{f^N}\|<R(\bk)$ (resp. $\rho(\Delta_{f^N})< R(\bk)$.)
Moreover, if $\Char \,\widetilde{\bk}=0$, we may choose $N=1$, if $\Char \,\widetilde{\bk}=p>0$, we may choose $N$ to be a power of $p.$
\end{pro}

\proof[Proof of Proposition \ref{prodeltaiteratesmall}]
We prove this lemma for $\|\Delta_f\|$, the proof for  $\rho(\Delta_f)$ is similar.

If $\Char \,\widetilde{\bk}=0$, this lemma is trivial. Now we assume that $\Char \,\widetilde{\bk}=p>0.$
Observe that
$$\Delta_{f^p}=(f^*)^p-\id=(\id+\Delta_f)^p-\id=\sum_{i=1}^{p-1}\binom{p}{i}\Delta_f^i+\Delta_f^p.$$
Since $p| \binom{p}{i}$ for $i=1,\dots,p-1$, we have 
$$\|\Delta_{f^p}\|\leq \max\{ p^{-1}\|\Delta_f\|,\|\Delta_f\|^p\}.$$
Observe that if $\|\Delta_f\|\geq R(\bk)$, we have $p^{-1}\|\Delta_f\|\leq \|\Delta_f\|^p$, hence $\|\Delta_{f^p}\|\leq \|\Delta_f\|^p.$

Since $\|\Delta_f\|<1$, there exists $s\geq 1$ such that $\|\Delta_f\|^{p^s}<R(\bk).$
Then there exists $t\in \{1,\dots,s\}$ such that $\|\Delta_{f^{p^t}}\|<R$. Otherwise, we have 
$$\|\Delta_{f^{p^s}}\|\leq \|\Delta_{f^{p^{s-1}}}\|^p\leq \dots\leq \|\Delta_f\|^{p^s}<R(\bk)$$ which is a contradiction.
\endproof

\subsection{Analytic DML}
The following result can be think as an analytic version of DML for $f$ closed to $\id$.
\begin{thm}\label{thmandmlsmall}
Assume that $\bk$ if of characteristic $0$. Let $X=\sM(A)$ be a $\bk$-affioid space. Let $f: X\to X$ be an automorphism with $\rho(\Delta_f)<1.$ Let $V$ be a Zariski closed subset of $X$. Then for $x\in X(\bk)$, the set $\{n\geq 0|\,\, f^n(x)\in V\}$ is a finite union of arithmetic progressions.
\end{thm}
\proof
By Proposition \ref{prodeltaiteratesmall}, after replacing $f$ by a suitable iterate, we may assume that $\rho(\Delta_f)<R(\bk).$
By Theorem \ref{thminteraffinoid}, there is a morphism $\Phi: \D\to X$ satisfying $\Phi(n)=f^n(x)$ for every $n\in \Z_{\geq 0}\subseteq \D(\bk)$.
Assume that $\{n\geq 0|\,\, f^n(x)\in V\}$ is infinite.
Since $\{n\geq 0|\,\, f^n(x)\in V\}\subseteq \Phi^{-1}(V)$ and $\{n\geq 0|\,\, f^n(x)\in V\}$ is Zariski closed in $\D$, then $\Phi^{-1}(V)=\D$. It follows that $\{n\geq 0|\,\, f^n(x)\in V\}=\Z_{\geq 0}$. This concludes the proof.
\endproof

The above proof comes from a step of the proof of \cite[Theorem 1.3]{Bell2010}. In the proof of \cite[Theorem 1.3]{Bell2010}, Bell, Ghioca and Tucker basically proved the case where $X$ is the $p$-adic polydisk. Our proof is just a reformulation of theirs in the language of Berkovich space. It is easy to proof \cite[Theorem 1.3]{Bell2010}(=(1) of Theorem \ref{thmdmlknownex}) from Theorem \ref{thmandmlsmall}.
Indeed, for a suitable identification $\C\simeq \C_p$ for some prime $p\gg 0$, after replacing $f$ by a suitable positive iterate, we may find an $f$-invariant polydisc $U$ containing $x$ such that $\rho(\Delta_{f|_U})<1$. Then we conclude the proof by Theorem \ref{thmandmlsmall}.

In \cite{Benedetto2010}, Benedetto, Ghioca, Kurlberg and Tucker gave an example which shows that DML does not hold for the map $f:\D^2\to \D^2$ over $\Q_p$ defined by $(x,y)\mapsto (x+1, py).$ As the above map is not an automorphism. We ask the following question.
\begin{que}[Analytic DML]\label{queanlydml}
Assume that $\bk$ is of characteristic $0$. Let $X=\sM(A)$ be a $\bk$-affioid space. Let $f: X\to X$ be an automorphism. Does DML hold for $f$?
	\end{que} 
If $\widetilde{\bk}$ is contained in the algebraic closure of a finite field, easy to show that the above question has a positive answer.
In this case, after replacing $f$ by a suitable positive iterate, we may assume that the reduction $\widetilde{x}$ of $x$ is fixed by the reduction $\widetilde{f}$ of $f.$ Not hard to show that there is an affinoid subdomain $Y$ of $X$ containing $x$ whose image under the reduction map is $\tilde{x}$ and it is $f$-invaraint. One can show that $\rho(\Delta_{f|_Y^N})<1$ for a suitable $N\geq 1.$
Then we conclude the proof by Theorem \ref{thmandmlsmall}.

\medskip

When $\bk$ is of positive characteristic, the $p$-adic interpolation lemma (=Theorem \ref{thminteraffinoid}) does not hold.
But it is also interesting to ask an analytic version of $p$-DML (c.f. \cite[Conjecture 13.2.0.1]{Bell2016}). 
\begin{que}\label{quepdmllocalan}
	Assume that $\bk$ is of characteristic $p>0$. Let $X=\sM(A)$ be a $\bk$-affioid space. Let $f: X\to X$ be an automorphism (with $\rho(\Delta_f)<1$). Does $p$-DML hold for $f$, i.e. up to a finite set, $\{n\geq 0|\,\, f^n(x)\in V\}$ is a finite union of arithmetic progressions along with finitely many 
	sets taking form $$\{\sum_{i=1}^mc_ip^{l_in_i}|\,\, n_i\in \Z_{\geq 0}, i=1,\dots,m\},$$
	where $m\in \Z_{>1}, k_i\in \Z_{\geq 0}, c_i\in \Q$?
	\end{que}
	
	One may expect that some consequences of the $p$-adic interpolation lemma (=Theorem \ref{thminteraffinoid}) hold in positive characteristic.
An example is \cite[Conjecture 1.1]{Xie2023}. Here we state it in a more general form.
\begin{con}\label{conlocalinforbit}
	Assume that $\bk$ is algebraically closed and of characteristic $0$.  Let $f: X\to X$ be an automorphism with $\rho(\Delta_f)<1$. 
	If $f$ is not of finite order,
	then there is a non-empty open subset $U$ of $X$ such that for every $\bk$-point in $U$, the orbit of $x$ is infinite.
\end{con}
Conjecture \ref{conlocalinforbit} holds in the characteristic $0$ case by Theorem \ref{thminteraffinoid}. This conjecture can be viewed as an analogy of \cite{Amerik2008}. If it has positive answer, it should be helpful for the Zariski dense orbit conjecture. 

\section{Questions}\label{question}
In this sections, we give some questions related to the DML problems. 

\subsection{Special cases of DML}
It seem that the DML conjecture is difficult in general. So it is natural to study it in some special cases.
Here we list some special cases which I feel in particular interesting.

\medskip

DML is open even for endomorphims of smooth projective surfaces. Among these cases, we feel that the following two cases are significant:
\begin{points}
\item[(1)] DML for endomorphisms $f:\P^2\to \P^2$ over $\C$ of degree at least $2$;
\item[(2)] DML for endomorphisms of $f_1\times f_2$ on $\P^1\times \P^1$ over $\C$ where $f_1,f_2$ are endomorphisms of $\P^1$ of the same degree $d\geq 2$.
\end{points}
In (2), if $f_1$ and $f_2$ has different degree, it is not hard to show that  DML holds using an argument in \cite[Section 9.1]{Xie2017a}.

\medskip

In \cite{Xie2017a} and Theorem \ref{thmendoplaneoverc}, the author proved the DML conjecture for endomorphisms of $\A^2_{\C}.$
We suspect that the techniques in this proof can be generalized to study more general endomorphisms of affine varieties.
The following two cases are in particular interesting:
\begin{points}
\item[(3)] DML for endomorphisms of $\A^l_{\C}$, $l\geq 3$;
\item[(4)] DML for endomorphisms of affine surfaces over $\C$.
\end{points}
Note that the dimension $l$ case of (3) implies the order $l$ case of
Conjecture \ref{connonlinearsml}.

\medskip

Few cases of DML are known for rational maps.  Two interesting open cases are as follows:
\begin{points}
\item[(5)] DML for birational self-maps of $\P^2_{\C}$;
\item[(6)] DML for skew-linear map on $\P^1\times \A^l$ taking form 
$$(x,y)\mapsto (g(x), A(x)y)$$ where $g$ is an automorphism of $\P^1$ and 
$A(x)$ is a matrix in $M_{l\times l}(k(x))$.
\end{points}
Recall that (6) is known when $\deg g\geq 2$ \cite{Ghioca2018}.
Combining with \cite{Ghioca2018}, Case (6) implies Conjecture \ref{connonconsml}.

\medskip

In the DML conjecture, the point $x$ is a closed point. It is natural to ask
\begin{points}
\item[(7)] DML for non-closed points i.e. $x\in X$ with $\dim\{x\}\geq 1.$
\end{points}
It is easy to see that the non-closed point case (7) is implied by the original DML conjecture. 

\medskip

For DML in positive characteristic, recall that Yang proved the following significant result:
\begin{thm}
\cite[Theorem 1.02]{Yang2023} Let $\bk$ be a complete non-archimedian valuation field with $\Char\, \bk=p>0.$
Let $\bk^{\circ}$ be its valuation ring and $\bk^{\circ\circ}$ be the maximal ideal.
Assume that $f:\P^N_{\bk}\to \P^N_{\bk}$  is totally inseparable  and it is a lift of Frobenius i.e. it
takes the following form:
$$f: [x_0:\dots :x_N]\mapsto [x_0^q+g_0(x_0^p,\dots,x_N^p):\dots :x_N^q+g_N(x_0^p,\dots,x_N^p)]$$
where $q$ is a power of $p$ and $g_i$ are homogenous polynomials of degree $q/p$ in $\bk^{\circ\circ}[x_0,\dots,x_N].$
\end{thm}
In this result, $f$ is asked to having two conditions:
\begin{points}
\item $f$ is totally inseparable i.e. the morphism $f^*\Omega_{\P^N_\bk}\to \Omega_{\P^N_\bk}$ is the zero map;
\item $f$ is a lift of Frobenius i.e. the reduction $\tilde{f}: \P^N_{\tilde{\bk}}\to \P^N_{\tilde{\bk}}$ of $f$ is the $q$-Frobenius map.
\end{points}
Note that the second condition still make sense when $\bk$ is of mixed characteristic.

It is interesting to ask whether DML holds when one of the above conditions hold.
\begin{que}\label{queimpyang}Does DML holds in the following cases:  
\begin{points}
\item[(8)] the case $f$ is  totally inseparable;
\item[(9)] the case $f$ is a lift of Frobenius.
\end{points}
\end{que}
In \cite{Xie2018}, the author proved serval results for lifts of Frobenius including DML for backward orbits \cite[Theorem 1.7]{Xie2018} .
\cite[Theorem 1.7]{Xie2018} is stated in mixed characteristic case, but its proof also works in pure positive characteristic case.

\subsection{DML for other maps}
The DML conjecture is for rational self-maps of algebraic varieties. It is interesting to ask the same question for other maps.
\begin{que}\label{queDMLkahler}Does DML holds for meromorphic self-maps of compact K\"ahler manifolds?
\end{que}
It is even challenging for automorphisms. This case is known for automorphisms of projective manifolds \cite{Bell2010}, but the proof in \cite{Bell2010} relies on the $p$-adic method, which seems hard to be applied in the K\"ahler case.

\medskip

One may ask similar question for proper meromorphic self-maps of proper Berkovich spaces. 
\begin{que}\label{queDMLberko}Does DML holds for meromorphic self-maps of proper Berkovich spaces?
\end{que}

For the non-proper case, we recall the following two questions asked in Section \ref{sectionautaff}.
\begin{que}[=Question \ref{queanlydml}]\label{queanlydmlr}
Assume that $\bk$ is of characteristic $0$. Let $X=\sM(A)$ be a $\bk$-affioid space. Let $f: X\to X$ be an automorphism. Does DML hold for $f$?
	\end{que}

\begin{que}[=Question \ref{quepdmllocalan}]
	Assume that $\bk$ is of characteristic $p>0$. Let $X=\sM(A)$ be a $\bk$-affioid space. Let $f: X\to X$ be an automorphism (with $\rho(\Delta_f)<1$). Does $p$-DML hold for $f$, i.e. up to a finite set, $\{n\geq 0|\,\, f^n(x)\in V\}$ is a finite union of arithmetic progressions along with finitely many 
	sets taking form $$\{\sum_{i=1}^mc_ip^{l_in_i}|\,\, n_i\in \Z_{\geq 0}, i=1,\dots,m\},$$
	where $m\in \Z_{>1}, k_i\in \Z_{\geq 0}, c_i\in \Q$?
	\end{que}

\medskip

Define \emph{piecewise algebraic self-maps} as follows:
Let $X=\sqcup_{j\in J} X_j$ be a decomposition of $X$ into finitely many constructible subsets.
Let $f_j$ be a morphism $X_j\to X$
We call such $f: X\to X$ a \emph{piecewise algebraic self-map} on $X$.
It is clear that $f$ is continuous for the constructible topology.
It is clear that every endomorphism of $X$ is a piecewise algebraic self-maps.

This class of maps seems quite interesting.
They are algebraic in natural, but not exactly algebraic. 
Here are some examples of piecewise algebraic self-maps not coming from algebraic maps.
\begin{exe}
Let $X=M_{N\times N}\simeq \A^{N^2}$ be the space of $N\times N$-matrix. 
Every matrix $A$ has a unique decomposition $A=D(A)+N(A)$ where $D(A)$ is diagonalizable $N(A)$ and nilpotent. 
There is a unique piecewise algebraic self-maps $f: X\to X$ such that $f(A)=D(A).$
Since $f\neq \id$ and it is the identity on a Zariski open subset of $X$, this map is not coming from rational self-maps.
\end{exe}

\begin{exe}
Let $X:=\A^2$. Let $f$ be the unique piecewise algebraic self-maps such that $f|_{\{x=2y\}}: (x,y)\mapsto (y, x-y)$ and $f|_{X\setminus\{x=2y\}}: (x,y)\mapsto (y, x+y).$
This map associates to the piecewise linear recurrence sequences $A_n, n\geq 0$ such that 
$A_{n+2}=A_n+A_{n+1}$ if $2A_n\neq A_{n-1}$; and $A_{n+2}=A_n-A_{n+1}$ if $2A_n= A_{n-1}$.
\end{exe}

\begin{que}\label{quepiecewadml}
Does DML hold for piecewise algebraic self-maps over $\C$?
	\end{que}

\subsection{Other questions}
The DML conjecture is not a full generalization of the original Mordell-Lang conjecture. In particular, it considers only the forward orbit but not the backward orbit.  In an informal seminar, S-W Zhang asked the following question to the author.
\begin{que}(=\cite[Question 1.3]{Xie2018})\label{quezhang} Let $X$ be a quasi-projective variety over $\C$ and $F:X\to X$ be a finite endomorphism. Let $x$ be a point in $X(\C)$. Denote by $O^-(x):=\cup^{\infty}_{i=0}F^{-i}(x)$ the backward orbit of $x$. Let $V$ be a positively dimensional irreducible subvariety of $X$. If $V\cap O^-(x)$ is Zariski dense in $V$, what can we say about $V$?
\end{que}

We note that if $V$ is preperiodic, then $V\cap O^-(x)$ is Zariski dense in $V$. As the dynamical Manin-Mumford conjecture, the converse is not true.  
We have counterexamples even when $F$ is a polarized\footnote{An endomorphism $F:X\to X$ on a projective variety is said to be  polarized if there exists an ample line bundle $L$ on $X$ satisfying $F^*L=L^{\otimes d}, d\geq 2$.} endomorphism. The following example is given by Ghioca, which is similar to \cite[Theorem 1.2]{Ghioca2011}.
\begin{exe}(=\cite[Example 1.4]{Xie2018}) Let $E$ be the elliptic curve over $\C$ defined by the lattice $\Z[i]\subseteq \C.$
Let $F_1$ be the endomorphism on $E$ defined by the multiplication by $10$ and $F_2$ be the endomorphism on $E$ defined by the multiplication by $6+8i$.
Set $X:=E\times E$, $F:=(F_1,F_2)$ on $X$. Since $|10|=|6+8i|$, $F$ is a polarized endomorphism on $X.$ 
Let $V$ be the diagonal in $X$ and $x$ be the origin. We may check that $V\cap O^-(x)$ is Zariski dense in $V$, but $V$ is not preperiodic. 
\end{exe}

As a special case of Question \ref{quezhang}, 
the author proposed the following conjecture in \cite{Xie2018}.
\begin{con}(=DML for coherence backward orbits)\label{condmlrev} Let $X$ be a quasi-projective variety over $\C$ and $F:X\to X$ be a finite endomorphism. Let $\{b_i\}_{i\geq 0}$ be a sequence of points in $X(\C)$ satisfying $f(b_i)=b_{i-1}$ for all $i\geq 1$. Let $V$ be a positively dimensional irreducible subvariety of $X$. If the $\{b_i\}_{i\geq 0}\cap V$ is Zariski dense in $V$, then $V$ is periodic under $F$.
\end{con}
This conjecture is true for lifts of Frobenius \cite[Theorem 1.7]{Xie2018}.

\medskip
Following the general principle of ``unlikely intersection problem" \cite{Zannier2012}, we propose the following question:
\begin{que}\label{quehdimdml}Let $X$ be a quasi-projective variaty over $\C,$ $f: X\to X$ be an endomorphism. Let $Z,V$ be irreducible subvarieties of $X$ with $\dim Z+\dim V<\dim X$. 
	\begin{points}
		\item Can we describe the set $\{n\geq 0|\,\, f^n(Z)\cap V\neq \emptyset\}$? For example, is it a union of arithmetic progressions?
		\item If $(\cap_{n\geq 0}f^n(Z))\cap V$ is Zariski dense, can we describe $V$?
		\end{points}
\end{que}
One may also ask the same question for rational self-maps. 
The part (i) can be think as a higher dimensional generalization of the original statement of DML and (ii) can be think as a higher dimensional generalization the geometric form of DML. Question \ref{quehdimdml} seems quite difficult. It will be also interesting to have a weak version of (i) like the weak DML.

\newpage
\bibliography{dd}
\end{document}